\documentclass[review,onefignum,onetabnum]{siamart190516}

\usepackage{amssymb,amsbsy,amsmath,amsfonts,amssymb,amscd}
\usepackage{mathrsfs}
\usepackage{algorithm,algpseudocode}
\usepackage{setspace}
\usepackage{float}
\usepackage{graphicx}
\usepackage{xcolor}
\usepackage[shortlabels]{enumitem}
\usepackage{subfig}
\usepackage{graphicx}
\usepackage{color}

\newtheorem{assump}{Assumption} 
\newtheorem{remark}{Remark}

\headers{Low rank integrator for the LBE}{Zhiyan Ding, Lukas Einkemmer and Qin Li}

\title{Dynamical low-rank integrator for the linear Boltzmann equation: error analysis in the diffusion limit}

\author{Zhiyan Ding\thanks{Department of Mathematics, University of Wisconsin-Madison, Madison, WI, 53706, USA (zding49@math.wisc.edu).}
\and%
Lukas Einkemmer\thanks{Department of Mathematics, University of Innsbruck, Innsbruck, Austria (lukas.einkemmer@uibk.ac.at)}
\and
Qin Li\thanks{Department of Mathematics and Wisconsin Institute for Discovery, University of Wisconsin-Madison, Madison, WI, 53706, USA (qinli@math.wisc.edu).}
}

\newcommand{\MSX}{\mathsf{X}}

\newcommand{\MSV}{\mathsf{V}}
\newcommand{\MSS}{\mathsf{S}}
\newcommand{\MSL}{\mathsf{L}}
\newcommand{\MSK}{\mathsf{K}}
\newcommand{\MSA}{\mathsf{A}}
\newcommand{\MSC}{\mathsf{C}}
\newcommand{\MSQ}{\mathsf{Q}}
\newcommand{\MSP}{\mathsf{P}}
\newcommand{\MSD}{\mathsf{D}}
\newcommand{\MSXi}{\mathsf{\Xi}}
\newcommand{\MSGamma}{\mathsf{\Gamma}}
\newcommand{\MSSigma}{\mathsf{\Sigma}}
\newcommand{\MSPi}{\mathsf{\Pi}}
\newcommand{\rd}{\,\mathrm{d}}

\newcommand{\Proj}{\mathcal P}
\newcommand{\n}{\mathrm{n}}

\newcommand{\revise}[1]{{\color{black}{#1}}}
\ifpdf
\hypersetup{
  pdftitle={Low rank integrator for linear Boltzmann},
  pdfauthor={Z. Ding, L. Einkemmer and Q. Li}
}
\fi

\begin{document}

\maketitle

% REQUIRED
\begin{abstract}
    Dynamical low-rank algorithms are a class of numerical methods that compute low-rank approximations of dynamical systems. This is accomplished by projecting the dynamics onto a low-dimensional manifold and writing the solution directly in terms of the low-rank factors. The approach has been successfully applied to many types of differential equations. Recently, efficient dynamical low-rank algorithms have been proposed in~\cite{E18,EL18} to treat kinetic equations, including the Vlasov--Poisson and the Boltzmann equation. There it was demonstrated that the methods are able to capture the low-rank structure of the solution and significantly reduce numerical cost, while often maintaining high accuracy. However, no numerical analysis is currently available.

    In this paper, we perform an error analysis for a dynamical low-rank algorithm applied to the multi-scale linear Boltzmann equation (a classical model in kinetic theory) to showcase the validity of the application of dynamical low-rank algorithms to kinetic theory. The equation, in its parabolic regime, is known to be rank one theoretically, and we will prove that the scheme can dynamically and automatically capture this low-rank structure. This work thus serves as the first mathematical error analysis for a dynamical low-rank approximation applied to a kinetic problem.
\end{abstract}

% REQUIRED

% REQUIRED

\maketitle

\section{Introduction}
Kinetic equations are a class of model equations used to describe the statistical behavior of a large number of particles that follow the same physical laws. They have been widely used in many aspects of physics and engineering. These equations, despite having different terms for various specific kinds of particles, share a similar structure: The dynamics is described in phase space and the solution is thus a distribution function $u(t,x,v)$ that counts the {density} of particles at a particular time $t$, location $x$, and velocity $v$. One main challenge for the computation of kinetic equation comes from the fact that the equation is supported on phase space instead of physical space, and thus the dimensionality of the problem typically doubles.

During the past decade, numerous methods have been proposed to numerically solve kinetic equations and many of them succeed in reducing computational cost without sacrificing accuracy. The main body of work concentrates on developing fast solvers for the collision operators and overcoming the stability requirement of the time-discretization \cite{SPT,JIN1996449}, through either finding fast methods to treat the transport term \cite{filbet2003,sircombe2009valis,crouseilles2011,einkemmer2015}, or implementing the resulting discrete system efficiently~\cite{grandgirard2016,einkemmer2018comparison}. However, reducing the complexity due to the high dimensionality is largely left unaddressed. The difficulty here is rather clear. To compute PDEs, one needs to sample a certain amount of discrete points per dimension, and for kinetic equations that have high dimensionality, the degrees of freedom is large, driving up the numerical cost.

The viewpoint of degrees of freedom merely depending on the dimensionality of the problem was challenged in recent years, see~\cite{NONNENMACHER20081346} and references therein. In these works, the traditional approach is abandoned and one explores the structure of the intrinsic low-dimensional manifold where the PDE solution lives on. The method is designed to follow the flow of the dynamics and to identify the important features in the evolution. Since the equation is projected onto a solution manifold of lower dimensionality, only the core information is preserved and redundant information is neglected. The numerical cost thus depends mainly on the intrinsic dimensionality of the dynamics, rather than the total number of discrete points. The method was first utilized to deal with coupled ODE systems \cite{Lubich2014} and then was applied to a range of problems for which its accuracy has been demonstrated.

The success in dealing with ODE systems inspired the generalization to PDEs. More recently, a class of efficient numerical schemes were proposed for computing kinetic equations, including the Vlasov--Poisson system~\cite{Kormanna,EL18}, the Vlasov-Maxwell system~ \cite{EOP_19} and the classical Boltzmann equation~\cite{E18} in both collisionless and strong-collisional regimes. In these experiments, it was observed that both the sophisticated Landau damping phenomenon for the VP system~\cite{ter1965,Mouhot2011}, and the incompressible Navier-Stokes limit for the Boltzmann equation are captured rather accurately at very low cost. Different equations may require slightly modified dynamical low-rank algorithms to suit the specific structures of the equation, but the general approach is quite similar. The main features in the evolution are preserved by following the flow of the equation projected onto the low dimensional solution manifold.

However, despite the strong intuition and the many promising numerical experiments, rigorous mathematical analysis is largely absent, especially in the PDE setting. There are two major difficulties, one from the theoretical level, and one is numerical. First, we often do not have results to show that the true solution is indeed approximately of low-rank. This point is more subtle than one might think. In fact, a Fourier expansion is also a low-rank approximation. But, a Fourier expansion uses a fixed set of basis functions, and the prediction of the rank is higher than the intrinsic rank of the solution. In addition, showing the low rank structure by performing a Fourier decomposition relies on assuming smoothness of the solution, which for kinetic equations is certainly problematic. It should be noted, however, that for elliptic equations some results can be obtained \cite{dahmen2016}.

The second difficulty is numerical: How to show the method can capture the low-rank structure of the solution? For the ODE case an analysis has been performed in \cite{Koch2007} and related works, but this analysis only applies to the non-stiff case and is thus not applicable to PDEs. The only mathematical analysis of a dynamical low-rank scheme in the PDE setting we are aware of is found in \cite{ostermann2018convergence}. This paper considers a special situation where the stiffness originates only from the linear differential operator whose corresponding flow is exactly computable within the low-rank formulation. This separation suggests a splitting procedure that decouples the stiff and non-stiff dynamics, which permits a convergence analysis of the resulting method. This is not the situation for kinetic equations, and thus the analysis does not apply.

In this paper, we will consider a multi-scale linear Boltzmann equation. The equation is equipped with a parameter called the Knudsen number, and when it is small, the equation has a strong collision operator, and is in the diffusive scaling. It is a well known result that when this happens, the solution reduces to a rank $1$ function, namely the solution can be written as a multiplication of two functions, one for $x$ and one for $v$ respectively. The first challenge mentioned above therefore turns out to be trivial in this case. In this particular paper, we aim at solving the second challenge, which is to prove the dynamical low-rank algorithm can indeed capture this low rank structure. More specifically, we will show, through utilizing the Hilbert expansion, that in the small Knudsen number limit, the method automatically captures the rank $1$ behavior of the solution. These results imply that compared to the traditional approaches that require discretization of the full phase space, the essential degrees of freedom only scale as the number of grid points in physical space.

We should mention that the results shown here only address the low-rank dynamics of the linear Boltzmann equation. It is a showcase for the validity of the algorithm on kinetic equations, and serves as a stepping stone for future studies on more complicated kinetic equations. We also stress that although the linear Boltzmann equation itself is linear, the studied dynamical low-rank algorithm is nonlinear and involved, making the proof highly non-trivial. We report the analytical result showing the performance of different time-integrators, discussing both advantages and disadvantages in the asymptotic regimes.

\subsection{Multi-scale linear Boltzmann equation} \label{sec:rte-diffusion-limit}
We first give a quick overview of the equation. In $d>1$ dimension, let $(t,x,v)\in\mathbb{R}^+\times\Omega_x\times\mathcal{S}^{d-1}$, the equation writes:
\begin{equation}\label{E1.1}
\partial_tu(t,x,v) = \mathcal{L}u = -\frac{v}{\epsilon}\cdot\nabla_xu(t,x,v)+\frac{\sigma(x)}{\epsilon^2}\left(\rho(x)-u(t,x,v)\right)\,,
\end{equation}
where $\sigma(x)$ is the scattering cross-section, $\Omega_x\subset\mathbb{R}^d$, $\epsilon$ is called the Knudsen number, which represents the ratio between the mean free path and the typical domain length, and $v$ is the angular variable. Since $v \in \mathcal{S}^{d-1}$, the unit sphere in $d$ dimensional space, all $v$ have the same amplitude. The density $\rho$ is defined as follows
\begin{equation}\label{rho}
\rho(x)=\langle u(t,x,v)\rangle_v=\frac{1}{\vert \mathcal{S}^{d-1} \vert}\int_{\mathcal{S}^{d-1}}u(t,x,v) \rd {v}\,,
\end{equation} 
where $\vert \mathcal{S}^{d-1} \vert$ is the volume of $\mathcal{S}^{d-1}$ and $\langle\cdot\rangle_v$ denotes the average with respect to $v\in \mathcal{S}^{d-1}$. For the sake of convenience we shorten the notation 
\[
\rd \mu_v=\frac{1}{|\mathcal{S}^{d-1}|}\rd v\,,
\]
where $\rd \mu_v$ is the normalized measure in $v\in\mathcal{S}^{d-1}$. The equation is written in diffusion scaling, meaning the transport term $v\cdot\nabla_xu$ and the collision term $\rho-u$ are enlarged by $\frac{1}{\epsilon}$ and $\frac{1}{\epsilon^2}$, respectively.

It is a well-known result that when $\epsilon \to 0$, the solution $u$ of \eqref{E1.1} will asymptotically approach $\rho$ and thus becomes $v$ independent. Furthermore, $\rho$ solves the heat equation, with $C_d$ being a constant that depends on $d$ only:
\begin{equation}\label{DL}
\partial_t\rho = \nabla_x\cdot\left(\frac{1}{C_d\sigma(x)}\nabla_x\rho\right)\,,\quad (t\,,x)\in\mathbb{R}^+\times\Omega_x\,.
\end{equation}
For this reason, this is termed the parabolic regime, and the limit is called the diffusion limit. For this paper, this limit is particularly interesting because it essentially shows that $u(t,x,v)\sim\rho(t,x)$ is approximately of rank $1$.

There are two difficulties in computing this equation. First, when the Knudsen number $\epsilon$ is small, both the transport term and the collision operator are stiff. A traditional numerical method would require small time stepsize that resolves the stiffness. Second, equation \eqref{E1.1} is posed in an $2d-1$ dimensional phase space. Suppose the equation has $N_x$ grid points per spatial direction and $N_v$ grid points per velocity direction, then $N_x^d N_v^{d-1}$ floating point numbers are needed. The large increase in the degrees of freedom is usually referred to as the \textit{curse of dimensionality} in the literature.

There are indeed techniques developed to enlarge the stability region to overcome the stiffness problem as mentioned in the first challenge. In the kinetic framework, it is usually referred to as asymptotic preserving~\cite{SAP,DEGOND20105630,LM_MM,BLM:08,DP:11,LP:14,DPJS,LARSEN1987283}. Also see reviews~\cite{HU2017103,SPL,DP:14}. However, the investigation into the second problem is mostly open. The belief that numerical cost mostly depends on the number of grid points and thus the degrees of freedom has been so firm, and essentially has not been challenged in the literature till the proposal of the dynamical low-rank approximation.

\subsection{Dynamical low rank approximation}
Dynamical low rank approximation is a systematic approach to tackle the curse of dimensionality for time-dependent problems. Believing that the number of grid points usually over-represent the necessary numerical information, the main aim of the approximation is to explore the rank structure of the solution manifold. Instead of seeking the solution on a fixed set of grid points, the method evolves the low rank representation, and characterizes the flow of the solution on a low dimensional manifold.

There are many ways to decompose the solution into its low rank presentation. For kinetic equations in particular, it is rather straightforward to separate the physical space and the velocity space:
\begin{equation}\label{eq:lr-expansion}
u(t,x,v) = \sum^r_{i,j=1}S_{i,j}X_i(t,x)V_j(t,v)\,,
\end{equation}
so that the low rank factors $X_i$ and $V_j$ depend only on either $x$ or $v$ respectively. In the formula, $r$ is the rank, and is typically significantly smaller than $\text{min}(N_x,N_v)$. For this presentation, the number of degrees of freedom is $O(r \left(N^d_x+N^{d-1}_v\right))$.

Historically, dynamical low-rank approximations have been considered extensively in high dimensional problems arising in quantum mechanics. Finding a low rank approximation there makes the computation tractable, see \cite{meyer90tmc,meyer09mqd} and~\cite{Lubich2008,lubich15tii,Conte2010} for a mathematical treatment. The application of the method in a general setting is studied in \cite{Koch2007,Koch2010,lubich13dab,arnold2014approximation}, where both the matrix and the general tensor formats are investigated. In the initial development, one significant disadvantage was observed: the methods are not robust with respect to over-approximation. Several strategies were proposed to address the issue, including both regularization, and the projector splitting integrator. The latter is viewed as a major improvement~\cite{Lubich2014} in enhancing the robustness with respect to small singular values~\cite{Kieri2016}. This approach was later extended to various tensor formats \cite{lubich15tii,lubich15tio,Haegeman2016,lubich2018tucker}, and is the approach to be utilized in this paper.

The application of dynamical low-rank approximation for kinetic equations is relatively recent. Numerical experiments have been performed on the Vlasov--Poisson \cite{EL18,EL18_cons}, the Vlasov--Maxwell \cite{EOP_19}, and the Boltzmann equation \cite{E18}, and the numerical evidence is very promising. However, theoretical justification have been lacking.

In particular, from a mathematical point of view we need to answer two questions
\begin{enumerate}
\item When do the solutions have low-rank structures?
\item \revise{Can dynamical low-rank approximation capture this structure?}
\end{enumerate}
The first question, in the kinetic framework, would be addressed by utilizing the fluid limit justification, as shown in~\eqref{E1.1} and~\eqref{DL}. Since one can show that as $\epsilon \to 0$ we get $u(t,x,v)=\rho(t,x)$, the velocity direction completely degenerates, and the rank of the representation in equation~\eqref{eq:lr-expansion} is simply $1$. It is then reasonable to expect that for small $\epsilon$ the solution is only slightly different from its rank-$1$ approximation.

To address the second question relies on proper algorithm-design. As mentioned above, we will use the projector-splitting approach. This leads to a set of three evolution equations for $S$, $X$, and $V$, respectively, and all three will be advanced in every single time step in a proper order. It turns out that the specific time integrator for each equation plays a crucial role. More specifically, we will show that the implicit Euler method, due to the lack of symmetry, requires fine time discretization for capturing the rank structure, while the~\revise{Crank--Nicolson--Implicit--Euler (CNIE)} method will converge with very relaxed constraints on the time step size but does require well-prepared initial data. \revise{Note that in the CNIE method we use the Crank--Nicolson scheme for $X$ and $S$ and the implicit Euler scheme for $V$.} We propose to use the implicit Euler method for the initial step and a very small time stepsize, to reveal the structure of the PDE solution, and then switch to~\revise{CNIE} to preserve this  structure.

The rest of the paper is organized as follows. In Section~\ref{sec:numlr} we present the application of the dynamical low-rank projector splitting integrator to the linear Boltzmann equation. Both the semi-discrete splitting scheme and the fully-discrete methods are presented, in Section~\ref{sec:semi_discrete} and~\ref{sec:full_discrete} respectively. In Section~\ref{sec:NA} we state and prove the main results. In Section~\ref{sec:cost} we study the numerical cost saving and an intuitive explanation of the error analysis is presented in Section~\ref{sec:intuition}. Section~\ref{sec:theorem1} and~\ref{sec:CN_analysis} are dedicated to the error analysis with the implicit Euler and~\revise{CNIE} time integrator respectively, and we summarize the analysis with a final proposal of an algorithm in Section~\ref{sec:3_summary}. Some parts of the proof in Section~\ref{sec:NA} are rather tedious, and we leave them to Appendix. Numerical evidences are presented in Section~\ref{sec:num}.

\section{Numerical scheme}\label{sec:numlr}
We present the application of the dynamical low-rank approximation method to the linear Boltzmann equation in this section. We follow the framework and notations in~\cite{EL18}. Throughout the paper we use the standard \revise{$L^2$} vector space in both the spatial and the velocity domain. That is, we have
\[
\langle f,g\rangle_x = \int f(x)g(x)\rd{x}\,,\quad \langle f,g\rangle_v = \int f(v)g(v)\rd\mu_v\,.
\]
We call a function $f(x,v)$ a rank-$r$ function if it can be expanded by a set of $r$ orthonormal basis functions in $x$ and $v$. All these functions are collected to the set $\mathcal{M}$:
\begin{definition}[Rank-$r$ function in $L_2(\rd x\rd v)$]
The collection of all rank-$r$ functions is denoted by
\begin{equation*}\label{M}
\mathcal{M}=\left\{f(x,v)\in L^2(\Omega_x\times\mathcal{S}^{d-1}): f(x,v) \text{ is rank-}r\right\},
\end{equation*}
where $f(x,v)$ is rank-$r$, meaning: $f(x,v) = \sum^r_{i,j=1}S_{i,j}X_i(x)V_j(v)$ where $X_i$ and $V_i$ being orthonormal in physical and velocity space respectively:
\begin{equation*}\label{E1.3}
\left\langle X_i\,,X_j\right\rangle_x=\int_{\Omega_x} X_iX_j \rd x=\delta_{ij},\quad \left\langle V_i\,,V_j\right\rangle_v=\int_{\mathcal{S}^{d-1}}V_iV_j\rd \mu_v =\delta_{ij}\,.
\end{equation*}
\end{definition}
We note that in this definition, only the rank, $r$, is fixed. The basis functions $X_i$ and $V_j$ can be arbitrary, as long as the orthogonality condition is satisfied. We further emphasize that $\mathcal{M}$ is not a function space; it is easily seen that the summation of two rank-$r$ functions may not be rank-$r$.

It is unlikely that the solution is of rank-$r$, i.e.~in $\mathcal{M}$, for all time. However, numerically one can argue that the solution is approximately of low rank. Thus for the numerical solution we seek a rank-$r$ approximation in $\mathcal{M}$ at every time step. The algorithm is consequently looking for a trajectory on the manifold $\mathcal{M}$ that resembles the evolution guided by the equation. In some sense, we need to project the equation in $L^2(\Omega\times\mathcal{S}^{d-1})$ to the manifold and find the equation that governs the dynamics of this trajectory in $\mathcal{M}$.

We denote by $u(t,x,v)$ the analytic solution and by $u_r(t,x,v)$ the numerical rank-$r$ approximation, then with the argument above, the governing equation for $u_r$ is
\begin{equation}\label{RD}
\partial_tu_r=\Proj_{u_r}(\mathcal{L}u_r) = \Proj_{u_r}\left(\frac{\sigma(x)}{\epsilon^2}\left(\rho-u_r\right)-\frac{v}{\epsilon}\cdot\nabla_xu_r\right)\,,
\end{equation}
where $\Proj_{u}f$ stands for the projection of $f$ onto the tangential plane of $\mathcal{M}$ at $u$. Denoting $\mathcal{T}_{u_r}\mathcal{M}$ the tangential plane of $\mathcal{M}$ at $u_r$, this projection operator ensures
\[
\dot{u}_r \in \mathcal{T}_{u_r}\mathcal{M}\,,
\]
guaranteeing that $u_r$ lies in the manifold $\mathcal{M}$ for all time.

We now look for explicitly expression of the tangential plane and the projection operator. For a nonlinear manifold $\mathcal{M}$, such expressions necessarily depend on the point where the tangent is looked for. Denote
\[
f=\sum_{ij}S_{ij}X_iV_j\in\text{Span}\{X_i\}_{i=1}^r\otimes\text{Span}\{V_j\}_{j=1}^r\,,
\]
then the tangential plane is given by
\begin{equation*}\label{TFM}
\begin{aligned}
&\mathcal{T}_f{\mathcal{M}}\\
&=\left\{g\in L^2(\rd x\rd v)\,:g=\sum^r_{i,j=1}X_i(x)\dot{S}_{i,j}V_j(v)+\dot{X}_i(x)S_{i,j}V_j(v)+X_i(x)S_{i,j}\dot{V}_j(v)\right.\\
&\left.\text{with}\ \dot{S}\in\mathbb{R}^{r\times r}, \dot{X}_i\in L^2(\Omega_x),\dot{V}_j\in L^2(\mathcal{S}^{d-1})\,,\text{and} \left\langle X_i,\dot{X}_j\right\rangle_x=\left\langle V_i,\dot{V}_j\right\rangle_v=0\right\}
\end{aligned}\,.
\end{equation*}
This set collects all functions whose infinitesimal, when added to $f$, still yields a rank-$r$ function.
In this definition, we notice that we are allowed to choose arbitrarily an $r\times r$ matrix $\dot{S}$, a function list $\dot{X}_i$ and a function list $\dot{V}_i$, as long as the gauge conditions, $\left\langle X_i,\dot{X}_j\right\rangle_x=\left\langle V_i,\dot{V}_j\right\rangle_v=0$, are satisfied. We note that the gauge conditions are imposed to guarantee the uniqueness of the low-rank factors. Interested readers are referred to~\cite{Lubich2014} for details.

With this definition, one has:
\begin{equation}\label{PF}
\Proj_fg=\Proj_{{X}}g+\Proj_{{V}}g-\Proj_{{V}}\Proj_{{X}}g,\quad \forall g\in L^2(\Omega_x\times\mathcal{S}^{d-1})\,,
\end{equation}
where the spatial and velocity projection are
\begin{align*}\label{PXgPVg}
\Proj_{{X}}g=\sum^r_i\langle X_i,g\rangle_x X_i\,,\quad \Proj_{{V}}g=\sum^r_i\langle V_i,g\rangle_v V_i\,,\quad \Proj_{{V}}\Proj_Xg=\sum^r_{i,j}X_i\langle X_iV_j,g\rangle_{x,v} V_j\,.
\end{align*}

Inserting~\eqref{PF} into~\eqref{RD} gives us the governing equation for $u_r$:
\begin{equation}\label{eqn:u_r}
\partial_tu_r=\Proj_{u_r}(\mathcal{L}u_r)=\left(\Proj_{X_{u_r}} + \Proj_{V_{u_r}} - \Proj_{V_{u_r}}\Proj_{X_{u_r}}\right)(\mathcal{L}u_r)\,.
\end{equation}
The numerical method will then be developed upon this formulation. We discuss the semi-discrete (in time) and the fully-discrete schemes (both implicit Euler and~\revise{CNIE}) in details in the following subsections.

%This will give us a formulation for one time step. To devise an asymptotic preserving scheme for integrating from $t=0$ to the final time $t_{\text{max}}$ we need a thorough understanding of the error analysis. In fact, the implicit Euler and Crank-Nicolson behave slightly \textcolor{black}{differently}. To optimally suit our purpose we run the Euler for one step before shifting to Crank-Nicolson scheme. We defer a detailed discussion to the end of Section 3.

\subsection{Semi-discrete low rank splitting method}\label{sec:semi_discrete}

In this section, we develop a projector-splitting method to solve equation \eqref{eqn:u_r}. In order to be concise, we simply denote the low rank solution $u_r$ by $u$ and we decompose the solution using its low-rank representation:
\begin{equation}\label{udecomposition}
\begin{aligned}
u(t,x,v)&=\sum^r_{i,j=1}X_i(t,x)S_{i,j}(t)V_j(t,v)\\
&=\mathrm{X}(t,x)\mathrm{S}(t)\mathrm{V}^\top(t,v)=\mathrm{X}(t,x)\mathrm{L}(t,v)=\mathrm{K}(t,x)\mathrm{V}^\top(t,v)
\end{aligned}\,,
\end{equation}
where $\mathrm{X}$ and $\mathrm{V}$ collects the basis functions
\begin{align*}
\begin{cases}
\mathrm{X}(t,x)=\left[X_1(t,x),X_2(t,x),\dots,X_r(t,x)\right]\,,\\
\mathrm{V}(t,v)=\left[V_1(t,v),V_2(t,v),\dots,V_r(t,v)\right]\,.
\end{cases}
\end{align*}
We will also be using quantities $\mathrm{K}$ and $\mathrm{L}$:
\begin{align}\label{msfL}
&\mathrm{L}(t,v)=\mathrm{S}(t)\mathrm{V}^\top(t,v)=\left[L_1(t,v),L_2(t,v),\dots,L_r(t,v)\right]^\top,\\
&\mathrm{K}(t,x)=\mathrm{X}(t,x)\mathrm{S}(t)=\left[K_1(t,x),K_2(t,x),\dots,K_r(t,x)\right]\,.
\end{align}
Computing~\eqref{eqn:u_r} at discrete times $t_n$ amounts to finding the governing equations that provide the updates of
\[
\mathrm{X}^n\,,\quad\mathrm{V}^n\,,\quad\text{and}\quad\mathrm{S}^n,
\]
respectively, for all $t_n$ with $n\geq 1$.

To specify the initial data, we project the initial condition onto $\mathcal{M}$ using the singular value decomposition (SVD), namely:
\begin{equation}\label{INITIAL}
    u(t=0,x,v)\approx \sum^r_{i,j=1} X_{i}(t=0,x)S_{i,j}(t=0)V_{j}(t=0,v) = \mathrm{X}^0(x)\mathrm{S}^0(\mathrm{V}^0(v))^\mathrm{T}\,.
\end{equation}

For updating $\mathrm{X,V,S}$, the Lie-Trotter splitting is used. From time step $t_n$ to $t_{n+1} = t_n+\Delta t$, we split the three operators on the right hand side of \eqref{eqn:u_r} into three sub-steps:
\begin{align}
\partial_tu&=\Proj_{\mathrm{X}_u}\left(\frac{\sigma(x)}{\epsilon^2}\left(\rho-u\right)-\frac{v}{\epsilon}\cdot\nabla_xu\right),\label{PX}\\
\partial_tu&=-\Proj_{\mathrm{V}_u}\Proj_{\mathrm{X}_u}\left(\frac{\sigma(x)}{\epsilon^2}\left(\rho-u\right)-\frac{v}{\epsilon}\cdot\nabla_xu\right),\label{PVPX}\\
\partial_tu&=\Proj_{\mathrm{V}_u}\left(\frac{\sigma(x)}{\epsilon^2}\left(\rho-u\right)-\frac{v}{\epsilon}\cdot\nabla_xu\right).\label{PV}
\end{align}
This splitting \eqref{PX}-\eqref{PV} takes place for each time step. That is, given the numerical solution at time $t_n$ we update the solution to $t_{n+1}$ by solving the three equations one after another. All equations are advanced for a full time step $\Delta t$. Different from directly solving \eqref{eqn:u_r}, each sub-step only changes one part of the decomposition. In particular, the first splitting step \eqref{PX} preserves $\mathrm{X}$, the last \eqref{PV} preserves $\mathrm{V}$, and the middle step \eqref{PVPX} merely updates $\mathrm{S}$. This allows us to update the three components separately. Below we detail the evolution of each sub-step:
\begin{enumerate}[-,topsep=0pt]
\item Updating \eqref{PX}:

Starting with $u^n(x,v)=\mathrm{X}^n(x)\mathrm{S}^n\left(\mathrm{V}^n(v)\right)^\top$, we run \eqref{PX} for a full time step $\Delta t$, and we denote the result by $u^{n+1/3}(x,v)$. The step preserves $\mathrm{X}$ and thus:
\begin{equation*}\label{Xequal1}
\mathrm{X}^{n+1/3}(x)=\mathrm{X}^{n}(x)\,.
\end{equation*}
To update $\mathrm{S}^{n+1/3}$ and $\mathrm{V}^{n+1/3}$, we plug the low rank formulation \eqref{udecomposition} into \eqref{PX} and obtain, for $1\leq i\leq r$:
\begin{equation*}
\begin{aligned}
\partial_tL_i=-\frac{1}{\epsilon}\sum^r_{j=1}\left\langle X^n_i,v\cdot\nabla X^n_j\right\rangle_xL_j+\frac{1}{\epsilon^2}\sum^r_{j=1}\left\langle X^n_i,\sigma X^n_j\right\rangle_x\left(\langle L_j\rangle_v-L_j\right)\,,
\end{aligned}
\end{equation*}
By using equation \eqref{msfL} we can simplify it to:
\begin{equation}\label{PXModify}
\partial_t\mathrm{L}(t,v)=-\frac{1}{\epsilon}\sum^d_{k=1}v_k\mathrm{A}^n_{\partial_k}\mathrm{L}(t,v)+\frac{1}{\epsilon^2}\mathrm{A}^n_{\sigma}\left(\left\langle\mathrm{L}(t,\cdot)\right\rangle_v-\mathrm{L}(t,v)\right),
\end{equation}
where $v=[v_1\,,\cdots\,, v_d]$, and $\mathrm{A}^n_{\partial_k}$ and $\mathrm{A}^n_{\sigma}$, both $\in\mathbb{R}^{r\times r}$, are matrix versions of the differential operator and the scattering operator respectively:
\begin{equation}\label{Asigma}
\begin{aligned}
&\left[\mathrm{A}^n_{\partial_k}\right]_{i,j}=\left\langle X^n_i\,,\partial_kX^n_j\right\rangle_x,\; (1\leq k\leq d)\,,\\
&\left[\mathrm{A}^n_\sigma\right]_{i,j}=\left\langle X^n_i,\sigma X^n_j\right\rangle_x\,.
\end{aligned}
\end{equation}
We denote the solution of~\eqref{PXModify} by $\mathrm{L}^{n+1/3}(v)$, and $\mathrm{S}^{n+1/3}$ and $\mathrm{V}^{n+1/3}(v)$ are obtained through the Gram-Schmidt process (QR factorization) which ensures the orthogonality of $\mathrm{V}^{n+1/3}(v)$:
\begin{equation*}\label{LQR}
\mathrm{L}^{n+1/3}(v)=\mathrm{S}^{n+1/3}\left(\mathrm{V}^{n+1/3}(v)\right)^\top\,.
\end{equation*}
Finally:
\begin{equation*}
u^{n+1/3}=(\mathrm{X}\mathrm{S}\mathrm{V}^\top)^{n+1/3}\,.
\end{equation*}

\item Updating \eqref{PVPX}:
In this step, equation \eqref{PVPX} is evolved for $\Delta t$ with initial condition $u^{n+1/3}$. Since the step only changes $\mathrm{S}$, we immediately obtain 
\begin{equation*}\label{Xequal2}
\mathrm{X}^{n+2/3}=\mathrm{X}^{n+1/3}\,,\quad\text{and}\quad\mathrm{V}^{n+2/3}=\mathrm{V}^{n+1/3}\,.
\end{equation*}
Plugging the low rank representation~\eqref{udecomposition} into~\eqref{PVPX}, we obtain, in matrix form:
%for $1\leq i,j\leq r$,
%\begin{equation*}\label{PVPXF}
%\begin{aligned}
%&\partial_tS_{i,j}(t)=\frac{1}{\epsilon}\sum^r_{p,q=1}\left\langle X^n_iV^{n+2/3}_j,\left(v\cdot\nabla_xX^n_p\right)V^{n+2/3}_q\right\rangle_{x,v}S_{p,q}(t)\\
%-&\frac{1}{\epsilon^2}\sum^r_{p,q=1}\left\langle X^n_i,\sigma(x)X^n_p\right\rangle_xS_{p,q}(t)\left(\langle V^{n+2/3}_j\rangle_v\langle V^{n+2/3}_q\rangle_v-\delta_{j,q}\right)\,.
%\end{aligned}
%\end{equation*}
\begin{equation}\label{PVPXModify}
\partial_t\mathrm{S}(t)=\frac{1}{\epsilon}\sum^d_{k=1}\mathrm{A}^n_{\partial_k}\mathrm{S}(t)\mathrm{\Xi}^{n+2/3}_{v_k}-\frac{1}{\epsilon^2}\mathrm{A}^n_{\sigma}\mathrm{S}(t)\mathrm{\Gamma}^{n+2/3}\,,
\end{equation}
where $\mathrm{\Xi}_{v_k}$ and $\mathrm{\Gamma}$, both $\in\mathbb{R}^{r\times r}$, are the matrix versions of the multiplication operator associated to $v_k$ ($1\leq k\leq d$) and the density term, respectively:
\begin{align}
\left[\mathrm{\Xi}_{v_k}\right]_{i,j}=\left\langle V_i(\cdot),(\cdot)_kV_j(\cdot)\right\rangle_v\,, \quad\text{and}\quad\left[\mathrm{\Gamma}\right]_{i,j}=\langle V_i(\cdot)\rangle_v\langle V_j(\cdot)\rangle_v-\delta_{i,j}\,.\label{density}
\end{align}
We denote the computed update of~\eqref{PVPXModify} by $\mathrm{S}^{n+2/3}$ and
\begin{equation*}
u^{n+2/3}=(\mathrm{X}\mathrm{S}\mathrm{V}^\top)^{n+2/3}\,.
\end{equation*}

\item Updating \eqref{PV}:

In this step, \eqref{PV} is evolved for a full time step $\Delta t$ with initial value $u^{n+2/3}$. This step preserves $\mathrm{V}$. Thus,
\[
\mathrm{V}^{n+1}(v) = \mathrm{V}^{n+2/3}(v)\quad\Rightarrow\quad \mathrm{\Xi}^{n+1}=\mathrm{\Xi}^{n+2/3}\,,\quad\text{and}\quad\mathrm{\Gamma}^{n+1}=\mathrm{\Gamma}^{n+2/3}\,.
\]
To update $\mathrm{X}^{n+1}$ and $\mathrm{S}^{n+1}$, the low rank formulation \eqref{udecomposition} is plugging into \eqref{PX}. In the matrix form we have:
%Then for all $1\leq j\leq r$ we have
%\begin{equation*}\label{PVLF}
%\partial_t K_j(x,t)+\frac{1}{\epsilon}\sum^r_{i=1}\left\langle V^{n+1}_j,\left(v\cdot \nabla_xK_i\right)V^{n+1}_i\right\rangle_v
%=\frac{1}{\epsilon^2}\sum^r_{i=1}K_i(x,t)\left(\langle V^{n+1}_i\rangle_v\langle V^{n+1}_j\rangle_v-\delta_{i,j}\right),
%\end{equation*}
\begin{equation}\label{PVModify}
\partial_t\mathrm{K}(t,x)+\frac{1}{\epsilon}\sum^d_{k=1}\partial_k\mathrm{K}(t,x)\mathrm{\Xi}^{n+1}_{v_k}=\frac{1}{\epsilon^2}\mathrm{K}(t,x)\mathrm{\Gamma}^{n+1}\,.
\end{equation}
%with (since $\mathrm{V}^{n+2/3}=\mathrm{V}^{n+1}$)
%\[
%\mathrm{\Xi}^{n+1}=\mathrm{\Xi}^{n+2/3}\,,\quad\text{and}\quad\mathrm{\Gamma}^{n+1}=\mathrm{\Gamma}^{n+2/3}\,.
%\]
Solving this equation we obtain $\mathrm{K}^{n+1}(v)$, and the orthogonality of $\mathrm{X}^{n+1}$ is ensured through the Gram-Schmidt process:
\begin{equation*}\label{KQR}
\mathrm{K}^{n+1}(x)=\mathrm{X}^{n+1}\mathrm{S}^{n+1}\,.
\end{equation*}
\end{enumerate}
With these three steps completed, one finally arrives at the numerical solution at time $t_{n+1}$
\[
u^{n+1}=\left(\mathrm{X}\mathrm{S}\mathrm{V}^\top\right)^{n+1}\,.
\]

\begin{remark} We describe the method without incorporating special boundary conditions. The problem is assumed to be a Cauchy problem with infinite boundary. In the numerical examples and the theorems that we prove below, we used periodic boundary condition to eliminate the possible complication induced by the boundary. If Dirichlet-type boundary condition is provided, the incoming data may drive the solution away from the low rank approximation, and a boundary layer that sees drastic changes in both $x$ and $v$ space is needed to damp the fluctuation. Recent analysis in this direction can be found in~\cite{Guo,LLS,LLS_general,LLS_geometry} and the references therein.
\end{remark}

\subsection{Fully-discrete low rank splitting method}\label{sec:full_discrete}
We discretize the aforementioned splitting method in time and phase space. We denote
\begin{equation}\label{definitionofV}
\mathcal{X}=\left\{x_1,x_2\dots,x_{N_x}\right\}\,,\quad\mathcal{V}=\left\{v_1,v_2\dots,v_{N_v}\right\}
\end{equation}
the sets of discrete points in $\Omega_x$ and $\mathcal{S}^{d-1}$. The discrete solution can then be represented as follows
\[
\mathbb{R}^{N_x\times N_v}\ni u^n=\MSX^n\MSS^n(\MSV^n)^\top\,,
\]
where $\MSS^n = \mathrm{S}^n$ and we abuse the notation by calling
\begin{equation}\label{DIS}
\MSX^n=\left[X^n_1\ X^n_2\ \dots\ X^n_{r}\right]\in\mathbb{R}^{N_x\times r}\,, \quad\text{and}\quad \MSV^n=\left[V^n_1\ V^n_2\ \dots\ V^n_{r}\right]\in\mathbb{R}^{N_v\times r}
\end{equation}
with $X^n_i$ and $V^n_i$ denoting the $i$-th mode evaluated at the grid points.

We now perform an implicit time integration of the three equations~\eqref{PXModify},~\eqref{PVPXModify}, and~\eqref{PVModify} to overcome the stability issue introduced by the stiffness.

\begin{itemize}
\item[--] For updating~\eqref{PXModify}, $\MSX^{n+1/3} = \MSX^n$. Applying the implicit Euler scheme gives
\begin{equation}\label{FDL}
\left\{
\begin{aligned}
&\frac{\MSL^{n+1/3}-\MSL^{n}}{\Delta t}+\frac{1}{\epsilon}\sum^{d}_{k=1}\MSA^n_{\partial_k}\MSL^{n+1/3}\MSPi_{v_k}=\frac{\MSA^n_\sigma}{\epsilon^2}\MSL^{n+1/3}\mathsf{C}\,,\\
&\text{QR decomposition: }\MSL^{n+1/3}=\MSS^{n+1/3}(\MSV^{n+1})^\top\,.
\end{aligned}
\right.
\end{equation}
or applying the Crank--Nicolson method~\revise{(first step of CNIE)} gives
\begin{equation}\label{FDLnew}
\left\{
\begin{aligned}
&\frac{\MSL^{n+1/3}-\MSL^{n}}{\Delta t}+\frac{1}{\epsilon}\sum^{d}_{k=1}\MSA^n_{\partial_k}\left(\frac{\MSL^{n+1/3}+\MSL^{n}}{2}\right)\MSPi_{v_k}=\frac{\MSA^n_\sigma}{\epsilon^2}\frac{\MSL^{n+1/3}+\MSL^{n}}{2}\MSC\,,\\
&\text{QR decomposition: }\MSL^{n+1/3}=\MSS^{n+1/3}(\MSV^{n+1/3})^\top\,.
\end{aligned}
\right.
\end{equation}
In the equation we have used
\begin{equation}\label{CXPi}
\begin{aligned}
&\MSA^n_{\partial_k}=(\MSX^n)^\top \mathsf{D}_k\MSX^n\,,\quad \MSA^n_\sigma=\MSX^\top\MSSigma\MSX\,,\\
&\mathsf{\Pi}_{v_k}=\text{diag}(\mathcal{V}_k)\,,\quad \mathsf{C}=\frac{1}{N_v}ee^\top-\mathsf{I}_{N_v}\,,
\end{aligned}
\end{equation}
with $e=\left(1,1,\dots,1\right)^\top$. $\MSSigma=\text{diag}(\sigma(\mathcal{X}))$ is an $N_x\times N_x$ matrix with evaluations of $\sigma$ at the grid points in $\mathcal{X}$ assigned as diagonal entries and $\mathsf{D}_k$ is the discrete approximation of $\partial_{x_k}$. The specific form of $\mathsf{D}_k$ depends on the spatial discretization. We have used the simple rectangle rule for the integration in $\mu_v$. This determines the form of $\mathsf{C}$. Other numerical integral rules could also be applied and the specific form of $\mathsf{C}$ will change accordingly. Obviously $\MSA^n_{\partial_k}$ and $\MSA^n_\sigma$ are the discrete versions of \eqref{Asigma}.
\item[--] For updating~\eqref{PVPXModify}, we note that $\MSX$ and $\MSV$ are preserved:
\[
\MSX^{n+2/3} = \MSX^{n+1/3}\,,\quad\text{and}\quad\MSV^{n+2/3} = \MSV^{n+1/3}\,.
\]
The direct application of the implicit Euler scheme gives%
\begin{equation}\label{FDS}
\frac{\MSS^{n+2/3}-\MSS^{n+1/3}}{\Delta t}-\frac{1}{\epsilon}\sum^d_{k=1}\MSA^n_{\partial_k}\MSS^{n+2/3}\mathsf{\Xi}^{n+2/3}_{v_k}=-\frac{\MSA^n_\sigma}{\epsilon^2}\MSS^{n+2/3}\mathsf{\Gamma}^{n+2/3}\,.
\end{equation}
Similarly, if Crank--Nicolson~\revise{(second step of CNIE)} is used , the scheme reads
\begin{equation}\label{FDSnew}
\begin{aligned}
&\frac{\MSS^{n+2/3}-\MSS^{n+1/3}}{\Delta t}-\frac{1}{\epsilon}\sum^d_{k=1}\MSA^n_{\partial_k}\frac{\MSS^{n+2/3}+\MSS^{n+1/3}}{2}\MSXi^{n+2/3}_{v_k}\\
=&-\frac{\MSA^n_\sigma}{\epsilon^2}\frac{\MSS^{n+2/3}+\MSS^{n+1/3}}{2}\MSGamma^{n+2/3}\,.
\end{aligned}
\end{equation}
Here $\MSXi^{n+2/3}_{v_k}$ and $\MSGamma^{n+2/3}$, both $r\times r$, are discrete versions of \eqref{density}:
\begin{equation}\label{CVB}
\MSXi^{n+2/3}_{v_k}=(\MSV^{n+2/3})^\top\mathsf{\Pi}_{v_k}\MSV^{n+2/3}\,,\quad \MSGamma^{n+2/3}=(\MSV^{n+2/3})^\top\MSC \MSV^{n+2/3}\,.
\end{equation}
\item[--] Finally for updating~\eqref{PVModify}, we note 
\[
\MSV^{n+1} = \MSV^{n+2/3}\,,\quad \MSGamma^{n+1}=\MSGamma^{n+2/3}\,,\quad \MSXi^{n+1}=\MSXi^{n+2/3}\,.
\]
Defining:
\begin{equation}\label{FDK2}
\MSK^{n+2/3}=\MSX^{n+2/3}\MSS^{n+2/3}\,,
\end{equation}
we apply the implicit Euler method for
\begin{equation}\label{FDK}
\left\{
\begin{aligned}
&\frac{\MSK^{n+1}-\MSK^{n+2/3}}{\Delta t}+\frac{1}{\epsilon}\sum^d_{k=1}\mathsf{D}_k\MSK^{n+1}\MSXi^{n+1}_{v_k}=\frac{\MSSigma}{\epsilon^2}\MSK^{n+1}\MSGamma^{n+1}\,,\\
&\text{QR decomposition:}\ \MSX^{n+1}\MSS^{n+1}=\MSK^{n+1}\,.
\end{aligned}
\right.
\end{equation}
The Crank--Nicolson method will not be used in the last step and thus we do not specify it.%
\end{itemize}
We finalize the update for
\[
u^{n+1}=\MSX^{n+1}\MSS^{n+1}(\MSV^{n+1})^\top\,.
\]

\revise{We stress that the choice of $\mathsf{D}_k$ is not straightforward. Indeed, in the $\epsilon\to0$ limit, a term $\mathsf{D}_k\Sigma^{-1}\mathsf{D}_k$ will appear as the discretization of the Laplacian term. To guarantee the self-adjointness on the discrete level, $\mathsf{D}_k$ needs to be specifically chosen. For more details see Remark~\ref{rmk:choice_D}.} 

We summarize the one-time step integrators in Algorithm \ref{alg:Euler} and \ref{alg:CN}.
\begin{algorithm}[htb]
\caption{\textbf{One-timestep implicit Euler integrator}}\label{alg:Euler}
\setstretch{0.1}
\begin{algorithmic}
\State \textbf{Input:}
\State 1. time step $\Delta t$;
\State 2. Input data:  $\MSX^n$, $\MSS^n$ and $\MSV^n$.
\State \textbf{Run:}

\begin{itemize}
\item Construct $\MSA^n_{\partial_k}, \MSA^n_{\sigma}$ by \eqref{CXPi}; Compute $\MSL^n=\MSS^n(\MSV^n)^\top$;

\item[--] Update $\MSL^{n+1/3}$ by solving \eqref{FDL} with \revise{$\Delta t$};

\item Perform QR decomposition to obtain $\MSS^{n+1/3},\MSV^{n+1}$ from \eqref{FDL}; Construct $\MSXi^{n+1}_{v_k}$ and $\MSGamma^{n+1}$ by \eqref{CVB};

\item[--] Update $\MSS^{n+2/3}$ with \revise{$\Delta t$} by solving \eqref{FDS};

\item Construct $\MSK^{n+2/3}=\MSX^n\MSS^{n+2/3}$ by \eqref{FDK2};

\item[--] Update $\MSK^{n+1}$ with \revise{$\Delta t$} by solving \eqref{FDK};

\item Perform QR decomposition to obtain $\MSX^{n+1},\MSS^{n+1}$ from \eqref{FDK}; 
\end{itemize}

\State \textbf{Output:} Numerical solution $\MSX^{n+1},\MSS^{n+1},\MSV^{n+1}$;

\end{algorithmic}
\end{algorithm}

\begin{algorithm}[htb]
\caption{\textbf{One-timestep \revise{Crank--Nicolson--Implicit--Euler (CNIE)} integrator}}\label{alg:CN}
\setstretch{0.1}
\begin{algorithmic}
\State \textbf{Input:}
\State 1. time step $\Delta t$;
\State 2. Input data: $\MSX^n$, $\MSS^n$ and $\MSV^n$.
\State \textbf{Run:}

\begin{itemize}
\item Construct $\MSA^n_{\partial_k}, \MSA^n_{\sigma}$ by \eqref{CXPi}; Compute $\MSL^n=\MSS^n(\MSV^n)^\top$;

\item[--] Update $\MSL^{n+1/3}$ by solving \eqref{FDLnew} with $\Delta t$;

\item Perform QR decomposition to obtain $\MSS^{n+1/3},\MSV^{n+1}$ from \eqref{FDLnew}; Construct $\MSXi^{n+1}_{v_k}$ and $\MSGamma^{n+1}$ by \eqref{CVB};

\item[--] Update $\MSS^{n+2/3}$ with $\Delta t$ by solving \eqref{FDSnew};

\item Construct $\MSK^{n+2/3}=\MSX^n\MSS^{n+2/3}$ by \eqref{FDK2};

\item[--]Update $\MSK^{n+1}$ with $\Delta t$ by solving \eqref{FDK};

\item Perform QR decomposition to obtain $\MSX^{n+1},\MSS^{n+1}$ from \eqref{FDK};
\end{itemize}

\State \textbf{Output:} Numerical solution $\MSX^{n+1},\MSS^{n+1},\MSV^{n+1}$;

\end{algorithmic}
\end{algorithm}

\section{Properties of the numerical scheme}\label{sec:NA}
We investigate the properties of the numerical method in this section. In particular we will discuss the computational complexity and prove that the method captures the rank-$1$ structure in the parabolic regime.

\subsection{Computational complexity}\label{sec:cost}
To analyze the computational complexity is rather straightforward. Denote $r$ the rank, $N_x$ and $N_v$ the cardinality of $\mathcal{X}$ and $\mathcal{V}$ respectively. The matrices in the updating formula, $\MSX^n$, $\MSS^n$, $\MSV^n$ are computed by solving \eqref{FDL}, \eqref{FDS}, \eqref{FDK}. The cost is summarized in Table~\ref{table:cost}:

\begin{table}[ht]
\begin{center}
\begin{tabular}{|l|c|}
\hline
Operations&floating point operations(flops)\\
\hline\hline
Calculation of $\MSA^n_{\partial_k}$,$\MSA^n_\sigma$ \eqref{CXPi}&\revise{$O(r^2N_x)$}\\
\hline
Preparation of $\MSL^n$&\revise{$O(r^2N_v)$}\\
\hline
Total cost of preparation& \revise{$O(r^2(N_v+N_x))$}\\
\hline\hline
Calculation of $\MSPi_{v_k}\otimes \mathsf{A}^n_{\partial_k}$, $\mathsf{C}\otimes\MSA^n_\sigma$ &$O(r^2N_v)$\\
 \hline
 Solving $\MSL^{n+1}$~\eqref{FDL}&$O(r^3N^3_v)$\\
 \hline
Total cost of updating $\MSL$ &$O(r^3N^3_v)$\\
 \hline\hline 
QR decomposition to obtain $\MSS^{n+1/3}$, $\MSV^{n+1/3}$&\revise{$O(r^2N_v)$}\\
 \hline
 Calculation of $\MSXi^{n+2/3}_{v_k}$ and $\mathsf{\Gamma}^{n+2/3}$~\eqref{CVB} &\revise{$O(r^2N_v)$}\\
\hline
Calculation of $\MSXi^{n+2/3}_{v_k}\otimes \mathsf{A}^n_{\partial_k}$, $\mathsf{\Gamma}^{n+2/3}\otimes\MSA^n_\sigma$ &$O(r^3)$ \\
\hline
Solving $\MSS^{n+2/3}$ \eqref{FDSnew}&$O(r^6)$ \\
\hline
Total cost of updating $\MSS$ &\revise{$O(r^6+r^2N_v)$}\\
\hline\hline
Preparation of $\MSK^{n+2/3}$ \eqref{FDK2}&\revise{$O(r^2N_x)$}\\
\hline
Calculation of $\MSXi^{n+2/3}_{v_k}\otimes \mathsf{D}_k$, $\mathsf{\Gamma}^{n+2/3}\otimes\mathsf{\Sigma}$&$O(r^2N_x)$ \\
\hline
Solving $\MSK^{n+1}$ \eqref{FDK}&$O(r^3N^3_x)$\\
\hline
Total cost of updating $\MSK$ & $O(r^3N^3_x)$\\
\hline\hline
QR decomposition to obtain $\MSX^{n+1}$, $\MSS^{n+1}$&\revise{$O(r^2N_x)$}\\
\hline\hline
Total cost of one time step &$O(r^3(r^3+N^3_v+N^3_x))$\\
\hline
\end{tabular}
\caption{Floating point operation required for the dynamical low-rank algorithm. We note that in the computation of $\MSL^{n+1}$, $\MSS^{n+2/3}$ and $\MSK^{n+1}$ we used the classical LU decomposition. \revise{These steps can be sped up in implementation with iterative algorithms such as Krylov iteration. We leave out the discussion for this part of computational saving.}}\label{table:cost}
\end{center}
\end{table}

% \begin{itemize}
% \item Preparation: Calculation of $\MSA^n_{\partial_k}$,$\MSA^n_\sigma$ needs $O(r^2N^2_x)$ floating point operations (flops).

% \item Update $\MSV$: Calculation of $\MSL^n$ needs $O(r^3N_v)$ flops.

% Calculation of $\Pi_{v_k}\otimes \mathsf{A}^n_{\partial_k}$, $\mathsf{C}\otimes\MSA^n_\sigma$ need $O(r^2N_v)$ flops.

% Solving $\MSL^{n+1}$ needs $O(r^3N^3_v)$ flops.

% Using QR decomposition to obtain $\MSS^{n+1/3}$, $\MSV^{n+1/3}$ needs $O(N^3_v)$.

% \item Update $\MSS$: Calculation of $\MSXi^{n+2/3}_{v_k}$ and $\mathsf{\Gamma}^{n+2/3}$ needs $O(r^2N^2_v)$

% Calculation of $\MSXi^{n+2/3}_{v_k}\otimes \mathsf{A}^n_{\partial_k}$, $\mathsf{\Gamma}^{n+2/3}\otimes\MSA^n_\sigma$ needs $O(r^3)$ flops.

% Solving $\MSS^{n+2/3}$ needs $O(r^6)$ flops.

% \item Update $\MSX$: Calculation of $\MSK^{n+2/3}$ needs needs $O(r^3N_x)$ flops.

% Calculation of $\MSXi^{n+2/3}_{v_k}\otimes \mathsf{D}_k$, $\mathsf{\Gamma}^{n+2/3}\otimes\mathsf{\Sigma}$ needs $O(r^2N_x)$ flops.

% Solving $\MSK^{n+1}$ needs $O(r^3N^3_x)$ flops.

% Using QR decomposition to obtain $\MSX^{n+1}$, $\MSS^{n+1}$ needs $O(N^3_x)$.

% \end{itemize}
Because $r =  O(1)\ll \text{min}\{N_x\,,N_v\}$, in conclusion, we need $O(N^3_x+N^3_v)$ flops per time step, instead of the $O(N^3_xN^3_v)$ as required by the traditional solvers.

\subsection{Intuition of the error analysis}\label{sec:intuition}
Before describing and proving our results in detail, in this section we first give a relatively vague justification on why the method works. As described in the introduction, there are two points we need to make:
\begin{itemize}
\item Why the true solution has an approximate low rank structure? This question is a rather fundamental, and is independent of the method chosen: the rank structure of the solution purely depends on the governing PDE we are studying here.
\item Why the method keeps track of the low rank structure? This question concerns the behavior of the specific method (dynamical low-rank approximation) we choose to use here.
\end{itemize}

To answer the first question, we merely cite the following result.
\begin{theorem}[Theorem 2 of~\cite{BardosSantosSentis:84}]\label{thm:diff}
Denote $u^\epsilon$ the solution to the linear Boltzmann equation, for $(t,x,v)\in\mathbb{R}^+\times\Omega_x\times\mathcal{S}^{d-1}$:
\begin{equation}\label{eqn:rte_thm}
\partial_t u^\epsilon= \mathcal{L}u^\epsilon = \frac{\sigma}{\epsilon^2}\mathcal{L}_0u^\epsilon+\frac{1}{\epsilon}\mathcal{L}_1u^\epsilon\,,
\end{equation}
where $\mathcal{L}_0u = \rho-u$ and $\mathcal{L}_1u=-v\cdot\nabla_xu$. Then in the zero limit of $\epsilon$, the solution converges to the solution of the diffusion equation $\partial_t\rho = \nabla_x\cdot\left(\frac{1}{d\sigma}\nabla_x\rho\right)$ in the sense that
\[
\|u^\epsilon(t,x,v)-\rho(t,x)\|_{L_2(\rd x\rd\mu)}\leq C \epsilon \,,
\]
where $C$ has no dependence on $\epsilon$.
\end{theorem}
We note that periodic boundary conditions are imposed in the original theorem to avoid complications that may come from boundary layers. Essentially this theorem states that in the zero limit of $\epsilon$, $u^\epsilon(t,x,v)$ loses its velocity dependence, and the dynamics will purely be reflected in the physical space. In some sense:
\[
u^\epsilon(t,x,v) = \rho(x) +\mathcal{O}(\epsilon)\,,
\]
can be seen as the rank-$1$ approximation with the threshold set at any value bigger than $\epsilon$.

\textcolor{red}{
    The proof for the theorem uses the Hilbert expansion. Formally, one writes $u^\epsilon = u_0+\epsilon u_1+\epsilon^2 u_2+\cdots$ and performs an asymptotic analysis. Showing Theorem~\ref{thm:diff} rigorously then amounts to bounding $u_2$ uniformly in $\epsilon$ (this exposition is rather lengthy and we omit it here).}

\textcolor{red}{To answer the second question amounts to proving that the dynamical low-rank approximation captures the rank structure of the solution. To do so, we recall
\begin{equation}\label{eqn:rte_proj}
\partial_t u_r= \mathcal{P}_{u_r}\mathcal{L}u_r=\mathcal{P}_{u_r} \left[\frac{1}{\epsilon^2}\mathcal{L}_0u_r+\frac{1}{\epsilon}\mathcal{L}_1u_r\right].
\end{equation}
Notice that the equation differs from its continuous counterpart~\eqref{eqn:rte_thm} by the projection $\mathcal{P}_{u_r}$ on the right hand side. The key to the analysis is two-folded:
\begin{itemize}
\item Showing that the projection operator $\mathcal{P}_{u_r}$ does not disturb the rank-$1$ structure of the solution;
\item Showing that upon the projection, the leading order of the numerical solution follows the correct dynamics.
\end{itemize}}
\textcolor{red}{To show the former, we mimic the analysis on the PDE level, and perform the asymptotic analysis for~\eqref{eqn:rte_proj} by setting $u_r = u_{r,0}+\epsilon u_{r,1}+\cdots$. Then in the leading order:
\[
\mathcal{O}\left(\frac{1}{\epsilon^2}\right):\quad\mathcal{P}_{u_r}\mathcal{L}_0u_{r,0} = 0\,.
\]
Note that if we suppose $u_r\in\text{Span}\{X_i\}\otimes\text{Span}\{V_j\}$, then with the definition of $\mathcal{L}_0$,% calling $u_r\in\text{Span}\{X_i\}\otimes\text{Span}\{V_j\}$, then
\[
\mathcal{L}_0u_r=\rho-u_r = \langle u_r\rangle_v-u_r\in\text{Span}\{X_i\}\otimes\text{Span}\{V_j\}\quad\Rightarrow\quad \mathcal{P}_{u_r}\mathcal{L}_0u_{r,0} = \mathcal{L}_0u_{r,0} = 0\,.
\]
This suggests that $\mathcal{P}_{u_r}$, at least in the leading order, does not induce  new rank structure, and thus one has $u_{r,0}\in\mathrm{Null}\mathcal{L}_0$, making $u_r$ approximately rank-$1$. To show that $u_{r,0}$ follows the right equilibrium flow requires more delicate derivations in the higher order expansions, and we leave the details to Section~\ref{sec:NA}.}

\subsection{Asymptotic analysis of the implicit Euler method} \label{sec:theorem1}

To unify the notation, throughout this section, we denote by $u^n$ a matrix of size \revise{$N_x\times N_v$} with $u^n_{ij}$ being the numerical solution at $t_n,x_i,v_j$. Denote $e=[1,\cdots, 1]^\top$ a column vector of length $N_v$, $e_\n = e/\sqrt{N_v}$ the normalized version, then $\rho^n=u^ne$ is a vector of length $N_x$ representing the discrete version of $\rho$ at time $t_n$ for spatial grid points $\mathcal{X}$. The hope is to show that density $\rho^n$ solves equation~\eqref{DL} in the limit $\epsilon\to 0$. 

Our first result concerns the behavior of the implicit Euler method in the asymptotic regime. The corresponding result, stated in Theorem \ref{Theorem1}, is shown under the following technical assumption.
\begin{assump}\label{Assump0}
For all $n$, there exists an orthogonal matrix $\MSQ^n\in\mathbb{R}^{r\times r}$ such that
\begin{equation}\label{AE1.111}
\begin{aligned}
&\MSV^{n,*}=\MSV^{n}\MSQ^n\\
=&\left(e_\n\,,\mathcal{V}_{1,\n}+\epsilon a_1\,,\mathcal{V}_{2,\n}+\epsilon a_2\,,\dots,\mathcal{V}_{d,\n}+\epsilon a_d\,,V^*_{d+2},\dots,V^*_r\right) + O(\epsilon^2)\,,
\end{aligned}
\end{equation}
where the components satisfy (for all $1\leq k\leq d$ and $i\geq d+2$)
\begin{align} \label{Aorth1}
e^{\top}_\n a_k=O(\epsilon),\ \left(\mathcal{V}_{k,\n}\right)^{\top}V^{n+1,*}_i=O(\epsilon),\ e^{\top}_\n V^{n+1,*}_i=O(\epsilon^2)\,.
\end{align}
\end{assump}
\revise{Here $\mathcal{V}_i = [(v_1)_i\,,(v_2)_i\,,\cdots (v_{N_v})_i]^\top$ collects the $i$-th dimension component of all coordinates in $\mathcal{V}$. Since $v\in\mathbb{R}^d$, the subindex of $\mathcal{V}_i$ changes from $1$ to $d$. We further denote the normalization by $\mathcal{V}_{i,\n} = \frac{\sqrt{d}}{\sqrt{N_v}}\mathcal{V}_i$.}
% Here $\mathcal{V}_i$ is the $i$-th index of the coordinate set $\mathcal{V}$ and $\mathcal{V}_{i,\n} = \frac{\sqrt{d}}{\sqrt{N_v}}\mathcal{V}_i$ is the normalized version.
\begin{remark}[Intuition of Assumption \ref{Assump0}] \revise{To understand this assumption we link it to the proof for the diffusion limit in the continuous setting. For proving Theorem~\ref{thm:diff}, it is expected that
\[
u^\epsilon = u_0+\epsilon u_1+\epsilon^2 u_2\,,
\]
where $u_0$ has no $v$ dependence and $u_1\sim v\cdot \nabla_xu_0$. This means the leading order of the equation should be homogeneous in $v$ and the next order linearly reflects $v$ in every dimension. This is presented in the definition of $\MSV^{n,\ast}$ where $e_\n$ and $\mathcal{V}_i$ present homogeneity and linearity in the leading and the next order respectively.}
\end{remark}
There are direct consequences of this assumption. Define $\alpha^{n}=(\MSV^{n})^{\top}e_\n$ and $\alpha^{n,*}=(\MSV^{n,*})^{\top}e_\n$, then
\begin{equation}\label{Aorth2}
\alpha^{n}=\MSQ\alpha^{n,\ast},\quad \alpha^{n,*}=(\MSV^{n,*})^{\top}e_\n=\left(1+O(\epsilon^4),O(\epsilon^2),\dots,O(\epsilon^2)\right)^{\top}\,.
\end{equation}
Similarly, let
\[
\MSXi^{n,\ast}_{v_k}=(\MSV^{n,*})^\top\Pi_{v_k}\MSV^n=\MSQ^\top\MSXi^{n}_{v_k}\MSQ\,,\quad\text{and}\quad \MSGamma^{n,\ast}=(\MSV^{n,*})^\top\MSC\MSV^n=\MSQ^\top\MSGamma^{n}\MSQ\,,
\]
we have for all $j\leq N_v$
\begin{equation}\label{CVSTAR}
\left(\MSXi^{n,\ast}_{v_k}\right)_{1,j}=\left(\MSXi^{n,\ast}_{v_k}\right)_{j,1}=\frac{1}{\sqrt{d}}\delta_{k+1,j}+O(\epsilon)\,,
\end{equation}
and some straightforward calculation gives:
%\ql{did you miss some $\top$ below? also why do the second and the third line matter?}~\zd{I think no $\top$ is missed. The first equation should be $\MSXi^{n}_{v_m}\MSXi^{n}_{v_n}$.}\ql{how did you get \eqref{CVSTAR2}?}~\zd{This comes from $\MSV^{n,\ast}\alpha^{n,\ast}=e_\n$ and $e_\n \MSC=0$. I put it to the end.}
\begin{align}
&(\alpha^{n})^\top\MSXi^{n}_{v_m}\MSXi^{n}_{v_n}\alpha^{n}=(\alpha^{n,\ast})^\top\MSXi^{n+1,\ast}_{v_m}\MSXi^{n+1,\ast}_{v_n}\alpha^{n,\ast}=\frac{1}{d}\delta_{m,n}+O(\epsilon^2)\,,\label{eqn:consequence_assump}\\
&(\alpha^{n})^\top\alpha^{n}=(\alpha^{n,\ast})^\top\alpha^{n,\ast}=1+O(\epsilon^4)\,,\label{eqn:consequence_assump1}\\
&\MSV^{n}\alpha^{n}=\MSV^{n,\ast}\alpha^{n,\ast}=e_\n+O(\epsilon^2)\,,\label{eqn:consequence_assump2}\\
&(\alpha^{n})^\top\MSGamma^{n}=(\alpha^{n,\ast})^\top\MSGamma^{n,\ast}\MSQ^\top=O(\epsilon^2)\,.\label{CVSTAR2}
\end{align}

This assumption essentially says there are components in $\MSV^{n+1}$ that can represent $e_\n$ and all $\mathcal{V}_{k,\n}$. It is a rather mild condition. Indeed, assuming the initial data has this form, then with small time stepsize, \revise{namely if $\Delta t\ll(\Delta x)^2$,} the condition holds true at all later times, as explained in detail in Appendix \ref{SectionAssump}.
%This assumption is mild as the conditions imposed on the $a_i$ simply enforce the orthogonality constraint on the low-rank factors. We furthermore have the flexibility to choose our initial low-rank factors in such a form that Assumption \ref{Assump0} is satisfied. It can then be shown that for sufficiently small time step sizes the condition remains true at all later times. We explain this in detail in Appendix \ref{SectionAssump}.

\begin{theorem}\label{Theorem1} We employ the implicit Euler method to compute $u^{n+1}$ from $u^n$, using equations~\eqref{FDL},~\eqref{FDS} and~\eqref{FDK}, then, \textcolor{black}{under Assumption \ref{Assump0}, there is a constant $C$ independent of $\epsilon$ so that:}
\[ \color{black}
\| u^{n+1} - \rho_0^{n+1}e^\top \|_2 \leq C \epsilon\,,
\]
where $\rho_0^{n+1}$ solves
\begin{equation}\label{Firstdiffusionlimit}
\frac{\rho^{n+1}_0-\rho_0^{n+2/3}}{\Delta t} = \frac{1}{d}\sum^d_{k=1}\mathsf{D}_{{k}}\left(\MSSigma^{-1}\mathsf{D}_{{k}}\rho^{n+1}_0\right)
\end{equation}
and $\|\rho^{n+2/3}_0-\rho^{n}_0\|_{2}= O\left(\frac{(\Delta t)^2}{(\Delta x)^4}\right)$. This means in the $\epsilon\to0$ limit, the numerical solution is approximately rank $1$, with the density solving the diffusion equation according to the implicit Euler method.
\end{theorem}
\begin{remark} \textbf{The bad.}
We note that the implicit low-rank integrator based on the implicit Euler scheme has a strong requirement on $\Delta t$. It needs to be extremely small due to the $O\left(\frac{(\Delta t)^2}{(\Delta x)^4}\right)$ error term. For example, even if we choose $\Delta t \propto (\Delta x)^4$ the error is only $\mathcal{O}(\Delta t)$ and we end up with a global error of order $O(1)$. This renders this numerical method impractical for time integration.
\end{remark}
\begin{remark} \textbf{The good.} The method is able to capture the correct low-rank structure. For $\epsilon \to 0$ we have $u^{n+1} = \rho^{n+1}_0 e^\top$ approximately, meaning the solution has rank $1$ and is a constant in $v$. This is precisely the analytic result derived for the diffusion limit, shown in Theorem~\ref{thm:diff}. 
\end{remark}
\begin{remark}\textbf{On the rigor of the proof} The proof we provide below for the theorem is formal in the sense that we use the Hilbert expansion without tracing the constant $C$'s dependence on $\epsilon$. This proof, however, can be made rigorous with C's dependence on $\epsilon$ explicitly removed. Essentially one follows the flow of the current structure of the proof, and gives a bound at the second order expansion, as is done in the continuous setting~\cite{BardosSantosSentis:84}. It is significantly more tedious but provides little intuition, so we omit it from the current paper.
%In Appendix~\ref{sec:appen_rigor} we provide the rigorous proof using the Crank-Nicolson method. For the proof to be rigorous, one needs a slightly higher requirement on the initial data: the data needs to have a compatible form that links the leading order and the $\epsilon$ order. This relation is seen in the original proof on the continuous setting in~\cite{BardosSantosSentis:84} as well. The proof is significantly more tedious and longer since one has to trace higher order expansion and give a more explicit bound.
\end{remark}

\begin{proof}
Throughout the proof, all quantities of interests will be expanded using the following ansatz, with the subindex standing for the level in the asymptotic expansion.
\begin{equation}\label{eqn:expand}
p = p_0+\epsilon p_1+\epsilon^2 p_2+\cdots\,.
\end{equation}
We now proceed by performing the three steps in the low-rank projector splitting integrator.

\textbf{Step 1:} This step preserves $\MSX$ and updates $\MSL$. Thus, we plug the asymptotic expansion of $\MSL$ into equation \eqref{FDL} and obtain
\begin{equation}\label{E1.10}
\begin{cases}
O(1/\epsilon^2): &\MSA^n_\sigma \MSL^{n+1/3}_0\MSC=0\,,\\
O(1/\epsilon):&\sum^{d}_{k=1}\MSA^n_{\partial_k}\MSL^{n+1/3}_0\MSPi_{v_k}=\MSA^n_\sigma \MSL^{n+1/3}_1\MSC\,,\\
O(1): &\frac{\MSL^{n+1/3}_0-\MSL^{n}_0}{\Delta t}+\sum^{d}_{k=1}\MSA^{n}_{\partial_k}\MSL^{n+1/3}_1\MSPi_{v_k}=\MSA^n_\sigma \MSL^{n+1/3}_2\MSC\,.
\end{cases}
\end{equation}
Since $\MSA^n_\sigma$ is invertible, the equation for the leading order implies that $\MSL^{n+1/3}_0$ lies in the null space of $\MSC^\top$. According to equation~\eqref{CXPi} $\MSC$ is symmetric with null space $\text{span}\{e_\n\}$. Immediately:
\begin{equation}\label{Lfirstorder}
\MSL^{n+1/3}_0=l^{n+1/3}_0e^\top_\n\,,
\end{equation}
where the $r\times 1$ vector $l^{n+1/3}_0$ is yet to be determined. Solving the $O(1/\epsilon)$ level equation for $\MSL^{n+1/3}_1$ we have:
\begin{equation}\label{E1.11}
\MSL^{n+1/3}_1=-\sum^d_{k=1}(\MSA^n_{\sigma})^{-1}\MSA^{n}_{\partial_k}l^{n+1/3}_0e^\top_\n\MSPi_{v_k}+l^{n+1/3}_1e^\top_\n.
\end{equation}
where the term $l^{n+1/3}_1 e^\top_\n$ is an undetermined component lying in the null space of $\MSC$. To close the system we consider the $O(1)$ equation, and multiply $e_\n$ on both sides:
\begin{equation*}
\frac{\MSL^{n+1/3}_0e_\n-\MSL^n_0e_\n}{\Delta t}+\sum^{d}_{k=1}\MSA^n_{\partial_k}\MSL^{n+1/3}_1\MSPi_{v_k}e_\n=0\,.
\end{equation*}
Plugging~\eqref{E1.11} into this equation and using
\begin{equation*}
e^\top_\n\MSPi_{v_{k_1}}\MSPi_{v_{k_2}}e_\n=\frac{1}{d}\delta_{k_1,k_2}\,,\quad e^\top_\n\MSPi_{v_{k}}e_\n=0\,,
\end{equation*}
we obtain
\begin{equation*}
\frac{l^{n+1/3}_0-\MSL^n_0e_\n}{\Delta t}
    -\frac{1}{d}\sum^{d}_{k=1}\MSA^n_{\partial_k} (\MSA^n_{\sigma})^{-1}\MSA^{n}_{\partial_k}l^{n+1/3}_0=0\,.
\end{equation*}
Considering $u^{n+1/3}_0 = \MSX^n\MSL^{n+1/3}_0 = \MSX^n l^{n+1/3}_0e^\top_\n$ and $\rho^{n+1/3}_0 = u^{n+1/3}_0e_\n/\sqrt{N_v} =\MSX^n l^{n+1/3}_0/\sqrt{N_v}$, we obtain (using $\MSX^\top\MSX=\mathsf{I}$):
\begin{equation}\label{FinalL}
\frac{\rho^{n+1/3}_0-\rho^{n}_0}{\Delta t}-\frac{1}{d}\sum^{d}_{k=1}\MSX^n\MSA^n_{\partial_k}(\MSA^n_{\sigma})^{-1}\MSA^{n}_{\partial_k}(\MSX^{n})^\top\rho^{n+1/3}_0=0\,.
\end{equation}

Perform the QR decomposition of $\MSL^{n+1/3}$ to obtain $\MSV^{n+1/3}$ and $\MSS^{n+1/3}$, and since $\MSV$ will not change in later steps, we have, also according to~\eqref{CVB}:
\[
\begin{aligned}
\MSV^{n+1}&=\MSV^{n+2/3}=\MSV^{n+1/3}\,,\\
\MSGamma^{n+1}&=\MSGamma^{n+2/3}=\MSGamma^{n+1/3}=(\MSV^\top\MSC\MSV)^{n+1/3}\,,\\
\MSXi^{n+1}_{v_k}&=\MSXi^{n+2/3}_{v_k}=\MSXi^{n+1/3}_{v_k}=(\MSV^\top\MSPi_{v_k}\MSV)^{n+1/3}\,.
\end{aligned}
\]
For convenience, we use the superscript $(\cdot)^{n+1}$ uniformly for $\MSV,\MSGamma,\MSXi$ in the following discussion. According to Assumption \ref{Assump0}, there is $\MSQ\in\mathbb{R}^{r\times r}$ such that~\eqref{AE1.111} holds true for $\MSV^{n+1,\ast}=\MSV\MSQ$, then~\eqref{eqn:consequence_assump} holds true at the new time step.
%\begin{enumerate}[$\bullet$,topsep=0pt]
%\item Using~\eqref{Aorth2}, 
%\begin{equation}\label{eqn:alphastar}
%\alpha^{n+1,*}=(\MSV^{n+1,*})^{\top}e_\n=\left(1+O(\epsilon^4),O(\epsilon^2),\dots,O(\epsilon^2)\right)^\top\,;
%\end{equation}
%
%\item Similarly, let $\MSXi^{n+1,\ast}_{v_k}=\MSQ^\top\MSXi^{n+1}_{v_k}\MSQ$ and $\MSGamma^{n+1,\ast}=(\MSV^{n+1,\ast})^\top \MSC\MSV^{n+1,\ast}$, we have:
%\begin{equation}\label{CVSTAR}
%\left(\MSXi^{n+1,\ast}_{v_k}\right)_{1,j}=\left(\MSXi^{n+1,\ast}_{v_k}\right)_{j,1}=\frac{1}{\sqrt{d}}\delta_{k+1,j}+O(\epsilon)\,,\quad 1\leq j\leq N_v
%\end{equation}
%and
%\[
%\alpha^{n+1}\MSGamma^{n+1}=\alpha^{n+1,\ast}\MSGamma^{n+1,\ast}\MSQ^\top=O(\epsilon^4)\,.
%\]
%\item As a consequence,
%\begin{equation}\label{eqn:consequence_assump}
%\begin{aligned}
%&(\alpha^{n+1})^\top\MSXi^{n+1}_{v_m}\MSXi^{n+1}_{v_n}\alpha^{n+1}=(\alpha^{n+1,\ast})^\top\MSXi^{n+1,\ast}_{v_m}\MSXi^{n+1,\ast}_{v_n}\alpha^{n+1,\ast}=\frac{1}{d}\delta_{m,n}+O(\epsilon^2)\,,\\
%&(\alpha^{n+1})^\top\alpha^{n+1}=(\alpha^{n+1,\ast})^\top\alpha^{n+1,\ast}=1+O(\epsilon^4)\,,\\
%&\MSV^{n+1}\alpha^{n+1}=\MSV^{n+1,\ast}\alpha^{n+1,\ast}=e_\n+O(\epsilon^2)\,.
%\end{aligned}
%\end{equation}
%\end{enumerate}
%
%
%
%
Without loss of generality, in the following part of the proof, we ignore high order terms of $\epsilon$.

\textbf{Step 2:} This step preserves $\MSX$ and $\MSV$ and updates $\MSS$. We plug the asymptotic expansion of $\MSS$ into equation \eqref{FDS} to obtain
\begin{equation}\label{EPS}
\begin{cases}
O(1/\epsilon^2):& \MSA^n_\sigma \MSS_0^{n+2/3}\MSGamma^{n+1}=0\,,\\
O(1/\epsilon):& -\sum^d_{k=1}\MSA^n_{\partial_k}\MSS^{n+2/3}_0\MSXi^{n+1}_{v_k}=-\MSA^n_\sigma \MSS^{n+2/3}_1\MSGamma^{n+1}\,,\\
O(1):& \frac{\MSS^{n+2/3}_0-\MSS^{n+1/3}_0}{\Delta t}-\sum^d_{k=1}\MSA^n_{\partial_k}\MSS^{n+2/3}_1\MSXi^{n+1}_{v_k}=-\MSA^0_\sigma \MSS^{n+2/3}_2\MSGamma^{n+1}\,.
\end{cases}
\end{equation}
Since $\MSA^n_\sigma$ is invertible, the leading order equation implies that $\MSS_0^{n+2/3}$ lies in the null space of $\MSGamma^{n+1}$. Thus, we have 
\begin{equation}\label{Sfirstorder}
\MSS_0^{n+2/3}=s_0^{n+2/3}(\alpha^{n+1})^\top,
\end{equation}
where $s^{n+2/3}_0\in\mathbb{R}^{r\times 1}$ is yet to be determined.

We now follow the same strategy as in step 1. That is, we plug the expression for $\MSS_0^{n+2/3}$ into the equation of order $\mathcal{O}(1)$ and then project out the null space of $\MSGamma^{n+1}$. 
\[
\MSS^{n+2/3}_1=-\sum^d_{k=1}(\MSA^n_\sigma)^{-1}\MSA^n_{\partial_k}\MSS^{n+2/3}_0\MSXi^{n+1}_{v_k}+s^{n+2/3}_1(\alpha^{n+1})^\top.
\]
Then we close the system by plugging the result for $\MSS_1^{n+2/3}$ into the equation of order $\mathcal{O}(\epsilon)$. To do so, we first multiply $\alpha^{n+1}$ on both sides of the equation and use \eqref{CVSTAR2} for
\[
\frac{\MSS^{n+2/3}_0\alpha^{n+1}-\MSS^{n+1/3}_0\alpha^{n+1}}{\Delta t}-\sum^d_{k=1}\MSA^n_{\partial_k}\MSS^{n+2/3}_1\MSXi^{n+1}_{v_k}\alpha^{n+1}=0\,.
\]
Noticing~\eqref{Sfirstorder} and~\eqref{eqn:consequence_assump} we obtain
\[
\frac{\MSS^{n+2/3}_0\alpha^{n+1}-\MSS^{n+1/3}_0\alpha^{n+1}}{\Delta t}+\frac{1}{d}\sum^d_{k=1}\MSA^n_{\partial_k}(\MSA^n_\sigma)^{-1}\MSA^n_{\partial_k}s^{n+2/3}_0=0\,.
\]
Using \eqref{eqn:consequence_assump1}, we obtain

%
%
%This yields, using~\eqref{eqn:consequence_assump},\eqref{eqn:consequence_assump1},and \eqref{CVSTAR2}:\ql{are you sure? I think you used~\eqref{CVSTAR2} and the first equation of~\eqref{eqn:consequence_assump}...also~\eqref{eqn:consequence_assump} has another $\alpha$ on the left. where did it go below? Pls BE PRECISE! Both referees said the paper is too technical and they didn't want to even try to read it -- it is mostly because our own notations are horrible. This is simple algebra derivation -- we should not have made it so difficult to read. Besides the proof below does not need the $\ast$ information at all -- why did we even both putting it in the proof -- move everything to the assumption as the consequence there is much neater.}~\zd{I write down the detail here. If we include this, it seems redundant. But this is not very obvious.
%\\
%1. Multiplying $\alpha^{n+1}$ on both sides, we use \eqref{CVSTAR2} to obtain
%\[
%\frac{\MSS^{n+2/3}_0\alpha^{n+1}-\MSS^{n+1/3}_0\alpha^{n+1}}{\Delta t}-\sum^d_{k=1}\MSA^n_{\partial_k}\MSS^{n+2/3}_1\MSXi^{n+1}_{v_k}\alpha^{n+1}=0
%\]
%2. Plugging $S^{n+2/3}_1$ into it, we use \eqref{eqn:consequence_assump} to obtain
%\[
%\frac{\MSS^{n+2/3}_0\alpha^{n+1}-\MSS^{n+1/3}_0\alpha^{n+1}}{\Delta t}+\frac{1}{d}\sum^d_{k=1}\MSA^n_{\partial_k}(\MSA^n_\sigma)^{-1}\MSA^n_{\partial_k}s^{n+2/3}_0=0
%\]
%3. Using $\MSS_0^{n+2/3}=s_0^{n+2/3}(\alpha^{n+1})^\top$ and \eqref{eqn:consequence_assump1}, we obtain \eqref{eqn:S0}.
%}
%
\begin{equation}\label{eqn:S0}
\frac{\MSS^{n+2/3}_0\alpha^{n+1}-\MSS^{n+1/3}_0\alpha^{n+1}}{\Delta t}+\frac{1}{d}\sum^{d}_{k=1}\MSA^n_{\partial_k}(\MSA^n_{\sigma})^{-1}\MSA^{n}_{\partial_k}\MSS^{n+2/3}_0\alpha^{n+1}=0\,.
\end{equation}
Since we further have $\MSV^{n+1}\alpha^{n+1}=e_\n$ by \eqref{eqn:consequence_assump2}, similar to \eqref{FinalL}, we finally obtain
\begin{equation}\label{FinalS}
\frac{\rho^{n+2/3}_0-\rho^{n+1/3}_0}{\Delta t}+\frac{1}{d}\sum^{d}_{k=1}\MSX^n\MSA^n_{\partial_k}(\MSA^n_{\sigma})^{-1}\MSA^{n}_{\partial_k}(\MSX^n)^\top\rho^{n+2/3}_0=0\,,
\end{equation}
where $\rho^{n+2/3}_0=\MSX^{n}_0\MSS^{n+2/3}_0(\MSV^{n+1})^\top e_\n/\sqrt{N_v}$. We then update the leading order of $\MSK$ according to%
\begin{equation}\label{K23}
\MSK^{n+2/3}_0=\MSX^{n}\MSS^{n+2/3}_0=\MSX^{n}s^{n+2/3}_0(\alpha^{n+1})^\top\,.
\end{equation}

\textbf{Step 3:} This step preserves $\MSV$ and updates $\MSK$. Thus, we plug the asymptotic expansion of $\MSK$ into equation \eqref{FDK} and obtain
\begin{equation}\label{eqn:FOK}
\begin{cases}
O(1/\epsilon^2):& \MSSigma\MSK^{n+1}_0\mathsf{\Gamma}^{n+1}=0\,,\\
O(1/\epsilon):&\sum^d_{k=1}\mathsf{D}_k\MSK^{n+1}_0\MSXi^{n+1}_{v_k}=\MSSigma\MSK^{n+1}_1\mathsf{\Gamma}^{n+1}\,,\\
O(1):& \frac{\MSK^{n+1}_0-\MSK^{n+2/3}_0}{\Delta t}+\sum^d_{k=1}\mathsf{D}_k\MSK^{n+1}_1\MSXi^{n+1}_{v_k}=\MSSigma\MSK^{n+1}_2\MSGamma^{n+1}\,.
\end{cases}
\end{equation}
Since $\MSSigma$ is invertible, the leading order equation implies that $\MSK^{n+1}_0$ lies in the null space of $\MSGamma^{n+1}$. Thus, we have 
\begin{equation}\label{K1}
\MSK^{n+1}_0=k^{n+1}_0(\alpha^{n+1})^\top\,,
\end{equation}
where $k^{n+1}_0\in\mathbb{R}^{N_x\times 1}$ is yet to be determined. We now follow the same strategy as in step 1 and 2. That is to plug the expression for $\MSK_0^{n+1}$ into the equation of order $\mathcal{O}(1)$ and then project out the null space of $\MSGamma^{n+1}$. Then we close the system by plugging the result for $\MSK_1^{n+2/3}$ into the equation of order $\mathcal{O}(\epsilon)$. This yields, once again using~\eqref{eqn:consequence_assump},\eqref{eqn:consequence_assump1},and \eqref{CVSTAR2}:
\begin{equation}\label{E1.13}
    \frac{\MSK^{n+1}_0\alpha^{n+1}-\MSK^{n+2/3}_0\alpha^{n+1}}{\Delta t}-\frac{1}{d}\sum^d_{k=1}\mathsf{D}_{k} \MSSigma^{-1} \mathsf{D}_k \MSK^{n+1}_0 \alpha^{n+1}=0.
\end{equation}
Now, since $\rho^{n+1}_0=u^{n+1}_0e_\n/\sqrt{N_v}$ and $u^{n+1}_0= \MSK_0^{n+1} \left(\MSV^{n+1}\right)^\top  = k^{n+1}_0\alpha^\top\left(\MSV^{n+1}\right)^\top$, we obtain, by using equation~\eqref{eqn:consequence_assump2},
\begin{equation}\label{AE1.21}
\frac{\rho^{n+1}_0-\rho_0^{n+2/3}}{\Delta t}-\frac{1}{d}\sum^d_{k=1}\mathsf{D}_{{k}}\left(\MSSigma^{-1}\mathsf{D}_{{k}}\rho^{n+1}_0\right)=0\,,
\end{equation}
which concludes the proof for~\eqref{Firstdiffusionlimit}. To show $\|\rho^{n+2/3}_0-\rho^{n}_0\|_{2}= O\left(\frac{(\Delta t)^2}{(\Delta x)^4}\right)$, we denote $\mathcal{L}=\sum^{d}_{k=1}\MSX^n\MSA^n_{\partial_k}(\MSA^n_{\sigma})^{-1}\MSA^{n}_{\partial_k}(\MSX^n)^\top$, then equation~\eqref{FinalL} and~\eqref{FinalS} can be written as
\begin{equation}\label{K0to23}
\rho^{n+2/3}_0=\left(\mathsf{I}+\frac{\Delta t}{d}\mathcal{L}\right)^{-1}\left(\mathsf{I}-\frac{\Delta t}{d}\mathcal{L}\right)^{-1}\rho^{n}_0=\left(\mathsf{I}-\frac{(\Delta t)^2}{d^2}\mathcal{L}^2\right)^{-1}\rho^{n}_0\,.
\end{equation}
By using $\left(\MSX^n\right)^\top\MSX^n=\mathsf{I}$ we can bound $\mathcal{L}$ as follows
\[
\|\mathcal{L}\|_{2}\leq \frac{C_d}{(\Delta x)^2\min(\sigma(x))}\,,
\]
with $C_d$ being a constant depending on $d$ only. Therefore,
\[
\|\rho^{n+2/3}_0-\rho^{n}_0\|_{2}= O\left(\frac{(\Delta t)^2}{(\Delta x)^4}\right)\quad\Rightarrow\quad \|\rho^{n+2/3}-\rho^{n}\|_{2}= O\left(\frac{(\Delta t)^2}{(\Delta x)^4}\right)+O(\epsilon)\,.
\]

Finally, since
\[
u^{n+1}_0 =  k^{n+1}_0\alpha^\top\left(\MSV^{n+1}\right)^\top=k_0^{n+1} \left( e^\top_\n \MSV^{n+1} \right) (\MSV^{n+1}_0)^\top\,,
\]
we follow that $e_\n$ lies in the span of $\MSL_0^{n+1/3}$ and thus in the span of $\MSV^{n+1}_0=\MSV^{n+1/3}_0$. We thus have $u^{n+1}_0 = k_0^{n+1} e^\top_n$, as desired.
\end{proof}

\subsection{Asymptotic analysis of the \revise{CNIE} scheme}\label{sec:CN_analysis}

As observed in the previous section, the implicit Euler has a strong requirement on $\Delta t$ for the low-rank algorithm to preserve the rank structure. In Theorem \ref{Theorem2} below, we will show that the symmetry in the \revise{first two steps of CNIE} will bring the following advantage: the method can preserve the rank structure independent of the time step size $\Delta t$. However, we do assume well-prepared initial data.
% \begin{assump}\label{Assump2}
% For any $n>0$, there exits an orthogonal matrix $Q^{n+1}\in\mathbb{R}^{r\times r}$ such that
% \begin{equation}\label{AEMnew}
% \begin{aligned}
% &\MSV^{n+1,*}=\MSV^{n+1}Q\\
% =&\left(e_\n\,,\mathcal{V}_{1,\n}+\epsilon a_1\,,\mathcal{V}_{2,\n}+\epsilon a_2\,\dots,\mathcal{V}_{d,\n}+\epsilon a_d\,,V^*_{d+2},\dots,V^*_r\right)+ O(\epsilon^2)\,,
% \end{aligned}
% \end{equation}
% where the components satisfy (for all $1\leq k\leq d$ and $i\geq d+2$)\ql{we may need to delete the alpha line}
% \begin{align}
% e_\n^{\top}a_k=O(\epsilon),\ (\mathcal{V}_k)^{\top}V^{n+1,*}_i=O(\epsilon),\ e_\n^{\top}V^{n+1,*}_i=O(\epsilon^2)\,.\label{Aorthnew1}
% \end{align}
% \end{assump}
% A direct consequence of this assumption is, defining $\alpha^{n+1,*}=(\MSV^{n+1,*})^{\top}e_\n$, then
% \begin{equation}
% \alpha^{n+1,*}=(\MSV^{n+1,*})^{\top}e_\n=\left(1+O(\epsilon^4),O(\epsilon^2),\dots,O(\epsilon^2)\right)^{\top}\label{Aorthnew2}.
% \end{equation}

% The assumption is mild and we refer to section \ref{sec:theorem1} as well as Appendix \ref{Theorem2detail} and \ref{SectionAssump} for more details.\ql{if the assumption is the same as the previous one, delete all things above}

\begin{theorem}\label{Theorem2}
Apply the \revise{CNIE} method to compute $u^{n+1}$ from $u^n$, using equation~\eqref{FDLnew},~\eqref{FDSnew} and~\eqref{FDK}. Assuming that $u^n = \rho^{n}_0e^\top + \mathcal{O}(\epsilon)$, and under Assumption \ref{Assump0}, there is a constant $C$ independent of $\epsilon$ so that
\begin{equation}\label{PreservU} \color{black}
\|u^{n+1}-\rho_0^{n+1}e^\top\|_2 \leq C \epsilon\,,
\end{equation}
where
\begin{equation*}
\frac{\rho^{n+1}_0-\rho_0^{n}}{\Delta t}-\frac{1}{d}\sum^d_{k=1}\mathsf{D}_{{k}}\left(\MSSigma^{-1}\mathsf{D}_{{k}}\rho^{n+1}_0\right)=0\,.
\end{equation*}
This means with well-prepared initial data, the limiting scheme is the implicit Euler method applied on the diffusion equation~\eqref{DL}.
\end{theorem}

\begin{remark} \textbf{The good.} It is immediate that Theorem~\ref{Theorem2} differs from Theorem~\ref{Theorem1} in that we preserve $\rho$ in the limit $\epsilon\to 0$ independent of the time step size. Note that both methods share the last step, which is responsible for propagating the diffusion equation.
\end{remark}
\begin{remark} \textbf{The bad.} The \revise{CNIE} based low-rank algorithm used here only ``preserves" the asymptotic limit, but is not able to drive the solution to the low-rank space. As stated in Theorem \ref{Theorem2} we add the assumption that the initial value already has the corresponding structure of rank $1$, up to an error of $\mathcal{O}(\epsilon)$, to ensure the initial data is well-prepared at each time step.
\end{remark}
\begin{remark}\label{rmk:choice_D} \textbf{The choice of $\mathsf{D}_k$} We note that $\mathsf{D}_k(\MSSigma^{-1}\mathsf{D}_k)$ may not be self-adjoint operator in general if $\mathsf{D}_k$ is not chosen properly. In fact, if $\mathsf{D}_k$ is chosen as upwind type, the discretization is a shifted diffusion by one grid point in space. If $\mathsf{D}_k$ is chosen as central-scheme type, the self-adjoint property can be preserved, at the sacrifice of staggered behavior in the numerical solution. This type of problem is rather typical. A common way in the literature to overcome it is by introducing the even-odd decomposition, the even part goes to the limit and the odd part diminishes. The two components are evolved with different type of fluxes~\cite{QL}. Here numerically we simply balance the two types with a correctly chosen weight so that central scheme plays the major role in the diffusion limit. We leave the possible extension of introducing even-odd decomposition to future research.
\end{remark}
\begin{proof}
Most of the proof is very similar to the case of the implicit Euler scheme detailed in Theorem \ref{Theorem1}. Thus, we will only highlight the main differences here and refer the reader to the appendix for a more thorough exposition.

The asymptotic expansion for the first step is given by
\revise{\begin{equation*}
\begin{cases}
    O(1/\epsilon): &\frac{1}{2}\MSA^n_\sigma \left( \MSL^{n+1/3}_0 + \MSL^{n}_0 \right) \MSC=0\,,\\
O(1):&\frac{1}{2}\sum^{d}_{k=1}\MSA^n_{\partial_k} \left( \MSL^{n+1/3}_0 + \MSL^{n}_0 \right) \MSPi_{v_k}
    =\frac{1}{2}\MSA^n_\sigma \left( \MSL^{n+1/3}_1 + \MSL^{n}_1 \right) \MSC\,,\\
O(\epsilon): &\frac{\MSL^{n+1/3}_0-\MSL^{n}_0}{\Delta t}+\frac{1}{2}\sum^{d}_{k=1}\MSA^{n}_{\partial_k} \left( \MSL^{n+1/3}_1 + \MSL^{n}_1 \right) \MSPi_{v_k}=\frac{1}{2}\MSA^n_\sigma \left( \MSL^{n+1/3}_2 + \MSL^{n}_2 \right) \MSC\,.
\end{cases}
\end{equation*}}
\revise{From the equation of order $\mathcal{O}(1/\epsilon)$, we only know that $\MSL^{n+1/3}_0+\MSL^{n}_0$} lies in the null space of $\MSC^\top$; we can \emph{not} make a similar claim for $\MSL^{n+1/3}_0$. However, from the assumption $u^n = \rho^{n}_0e^\top + \mathcal{O}(\epsilon)$ and from $u^{n} = \MSX^{n} \MSL^{n}$ we immediately obtain $\MSL_0^n = l^{n}_0 e^\top_\n$.  This is an important ingredient in the remainder of the proof. It is also the first major difference between the present proof for the \revise{CNIE} and the proof of Theorem \ref{Theorem1}, where this condition is automatically satisfied independent of the chosen initial value. Using similar arguments as in the proof of Theorem \ref{Theorem1} we can then show (see Lemma \ref{lemmanew} in the appendix for details):
\begin{itemize}
\item $\MSL^{n+1/3}$ in \eqref{FDLnew} has the same form as \eqref{Lfirstorder}.
\item $\MSS^{n+1/3}$ in \eqref{FDSnew} has the same form as \eqref{Sfirstorder}.
\item The equations for computing $u_0^{n+1/3}$ and $u_0^{n+2/3}$ (equations \eqref{AnewFinalL} and \eqref{AnewFinalS} in the appendix).
\end{itemize}
With further details given in Appendix \ref{Theorem2detail}, we obtain
\begin{align*}
&\rho^{n+1/3}_0=\left(\mathsf{I}-\frac{(\Delta t)}{2d}\mathcal{L}^n\right)^{-1}\left(\mathsf{I}+\frac{(\Delta t)}{2d}\mathcal{L}^n\right)\rho^n_0,\\
&\rho^{n+2/3}_0=\left(\mathsf{I}+\frac{(\Delta t)}{2d}\mathcal{L}^n\right)^{-1}\left(\mathsf{I}-\frac{(\Delta t)}{2d}\mathcal{L}^n\right)\rho^{n+1/3}_0,\nonumber
\end{align*}
where as before $\mathcal{L}=\sum^{d}_{k=1}\MSX^n\MSA^n_{\partial_k}(\MSA^n_{\sigma})^{-1}\MSA^{n}_{\partial_k}(\MSX^n)^\top$. Combining these two equations we get $\rho^{n+2/3}_0=\rho^{n}_0$ for all $n>0$, which is the desired result.
\end{proof}

\subsection{Summary}\label{sec:3_summary}
Viewing Theorem~\ref{Theorem1} and Theorem~\ref{Theorem2}, it is clear that the implicit Euler method is able to capture the low-rank structure of the diffusion limit from arbitrary data at the cost of requiring a very small time step size, while the \revise{CNIE} scheme preserves the rank structure of the well-prepared initial data without constraint on the time step size.

When we design algorithms, we should make use of the benefit of both. To do so, one can simply run the implicit Euler scheme with a $\Delta t_1$ that is small enough such that the error term $O\left((\Delta t_1)^2/(\Delta x)^4\right)$ is controlled to the desired accuracy. This time integrator, however, gets applied only once, and then in all subsequent time steps, i.e.~all steps from time $t_1$ to time $t_\text{max}$, we apply the \revise{CNIE} time integrator. Theorem \ref{Theorem1} suggests that with the first step taken, all the prerequisites of Theorem \ref{Theorem2} are already guaranteed, so the rank structure will be preserved along the propagation, as Theorem~\ref{Theorem2} suggests. In conclusion, for $\epsilon \to 0$, the numerical solution is rank $1$ and homogeneous in velocity, with the density $\rho$ satisfying the implicit Euler discretization of the diffusion equation~\eqref{DL}.

We summarize the final method in Algorithm \ref{Algorithm1}.

\begin{algorithm}[htb]
    \caption{\textbf{(Dynamical low rank splitting method for the linear Boltzmann equation).}}
\label{Algorithm1}
\setstretch{0.1}
\begin{algorithmic}
\State \textbf{Preparation:}
\State 1. Initial data: $u(t=0,x,v)$.
\State 2. Input data:  final time: $t_\text{max}$; rank number: $r$; time step: $\Delta t_1\leq\Delta t_2$.
\State 3. Discretization points: $\mathcal{X}=\left\{x_1,x_2\dots,x_{N_x}\right\}, \mathcal{V}=\left\{v_1,v_2\dots,v_{N_v}\right\}$
\State 4. Initialization: use SVD to construct initial $u^0$ such that~\eqref{INITIAL}:
\[u^0=\MSX^0\MSS^0(\MSV^0)^{\top}\approx u(t=0,x,v)\]
\State 5. Construct $\MSPi_{v_k}$, $\MSC$ by \eqref{CXPi}.
\State \textbf{Run:} set $n=0$;
\State $\qquad$ call \textbf{One-timestep implicit Euler integrator} using $\Delta t_1$ to update $\MSX^1,\MSS^1,\MSV^1$, using Algorithm~\ref{alg:Euler}
\State \textbf{While} $t<t_\text{max}$: $n\to n+1$;
\State $\qquad$ call \textbf{One-timestep \revise{CNIE} integrator} using $\Delta t_2$ to obtain $\MSX^{n+1},\MSS^{n+1},\MSV^{n+1}$ using Algorithm~\ref{alg:CN}.
\State \textbf{end}
\State \textbf{Output:} Numerical solution $\MSX^n$, $\MSS^n$ and $\MSV^n$ for $t_n\leq t_\text{max}$.
\end{algorithmic}
\end{algorithm}

\section{Numerical experiments}\label{sec:num}

In this section we present some numerical evidence to showcase the property of the proposed dynamical low-rank integrator. Strictly speaking the linear Boltzmann equation is only defined for $d\geq 2$. However, in the quasi-2D case, assuming the data is homogeneous in $y$-direction, the equation degenerates to a problem with 1D in space. In this case the diffusion limit becomes:
\begin{equation}\label{DIFF1d}
  \partial_tu = \frac{1}{3}\partial_x\left(\sigma^{-1}(x)\partial_xu\right),\quad (t,x)\in\mathbb{R}^+\times\Omega_x\,.
\end{equation}
To measure the error we define:
\begin{align}\label{eqn:error_def_num}
\mathcal{P}^{error}_{x,l}=\|u-\mathsf{X}_l(\mathsf{X}_l)^\top u\|_F \quad\text{and}\quad \mathcal{P}^{error}_{v,l}=\|u-u\mathsf{V}_l(\mathsf{V}_l)^\top\|_F\,,
\end{align}
where $u$ is the reference solution and $\mathsf{X}$ and $\mathsf{V}$ are the numerical solution computed using the dynamical low-rank integrator. The error is measured in the Frobenius norm for the matrix, which is equivalent to $L_2(\rd x\rd\mu_v)$ in the continuous version.

\subsection{Projection error and singular value test}
\begin{figure}
     \centering
     \subfloat{\includegraphics[width=0.47\textwidth]{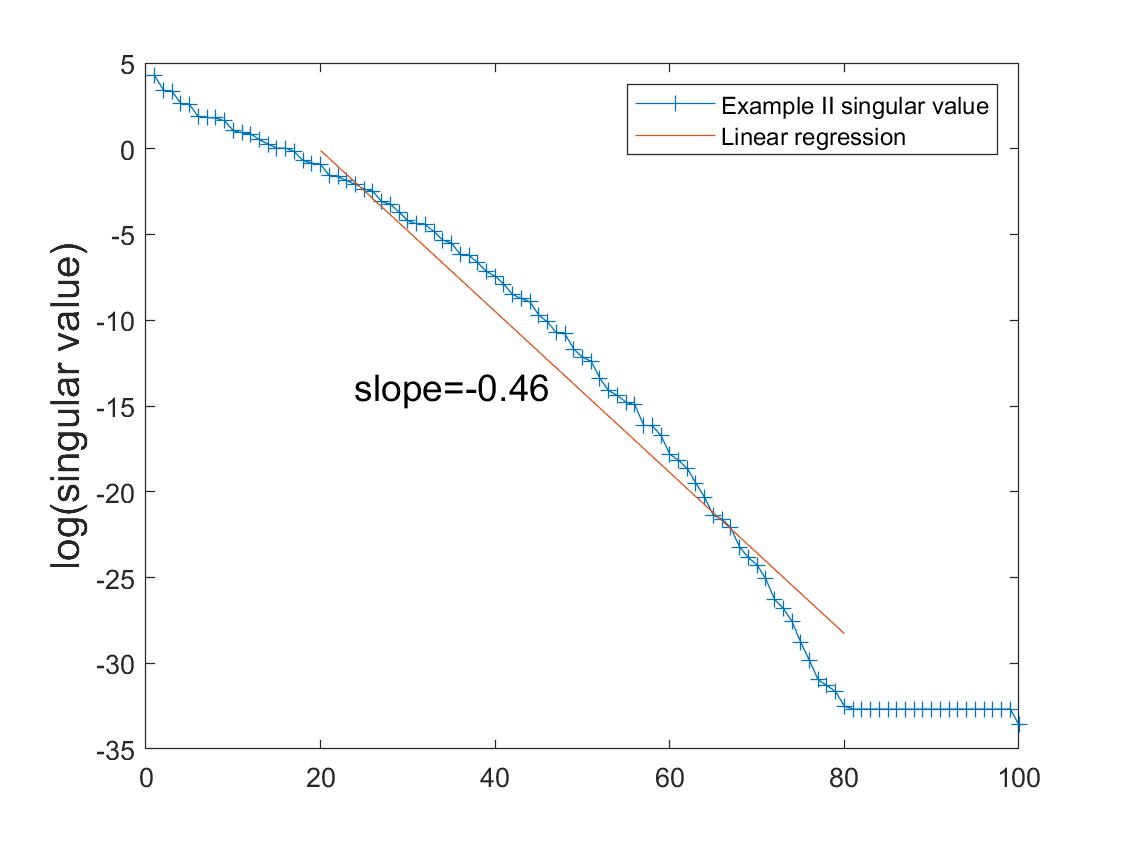}}
     \subfloat{\includegraphics[width=0.47\textwidth]{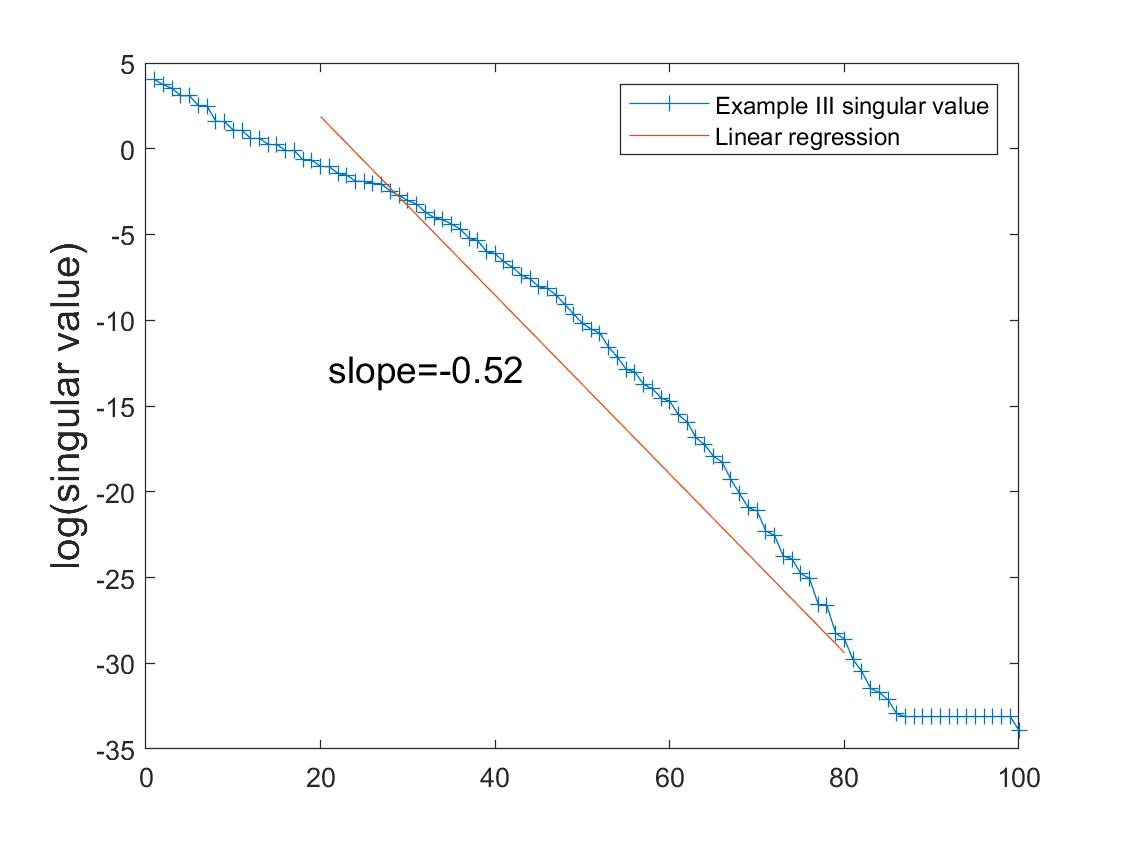}}
     \caption{Here we use \revise{$N_x=200, N_v=100, \Delta t=\Delta x/3$} and the direct implicit solver with upwind discretization for $\mathsf{D}$ to compute $u$. Left: Singular values of the solution $u$ in Example $\mathrm{II}$. Right: Singular values of the solution $u$ in Example $\mathrm{III}$.}
     \label{Figure7}
\end{figure}

Before testing the dynamical low-rank numerical integrator, we first numerically justify that the solution is indeed of low rank. Setting the initial data to
\begin{equation*}
f(0,x,v)=\left\{
	\begin{aligned}
     &\ 2,\quad 0.8<x<1.2\\
     &\ 0,\quad \text{otherwise},
	\end{aligned}
	 \right.
\end{equation*}
with the equation equipped with isotropic scattering and a varying cross sections (see~\eqref{eqn:sigma_example_2} and~\eqref{eqn:sigma_example_3} below) at $\epsilon=1$, we compute the equation with fine grids (\revise{$N_x=200$, $N_v=100$, and $\Delta t=\Delta x/3$}) till $t_{\max}=1$. We plot the singular values of the solution in log-scale in Figure~\ref{Figure7}. It is clear that the singular values decay exponentially fast for both cases. This gives us the foundation to believe that the dynamical low-rank approximation would work.

\subsection{Example $\mathrm{I}$} 
\begin{figure}
     \centering
     \subfloat{\includegraphics[width=0.47\textwidth]{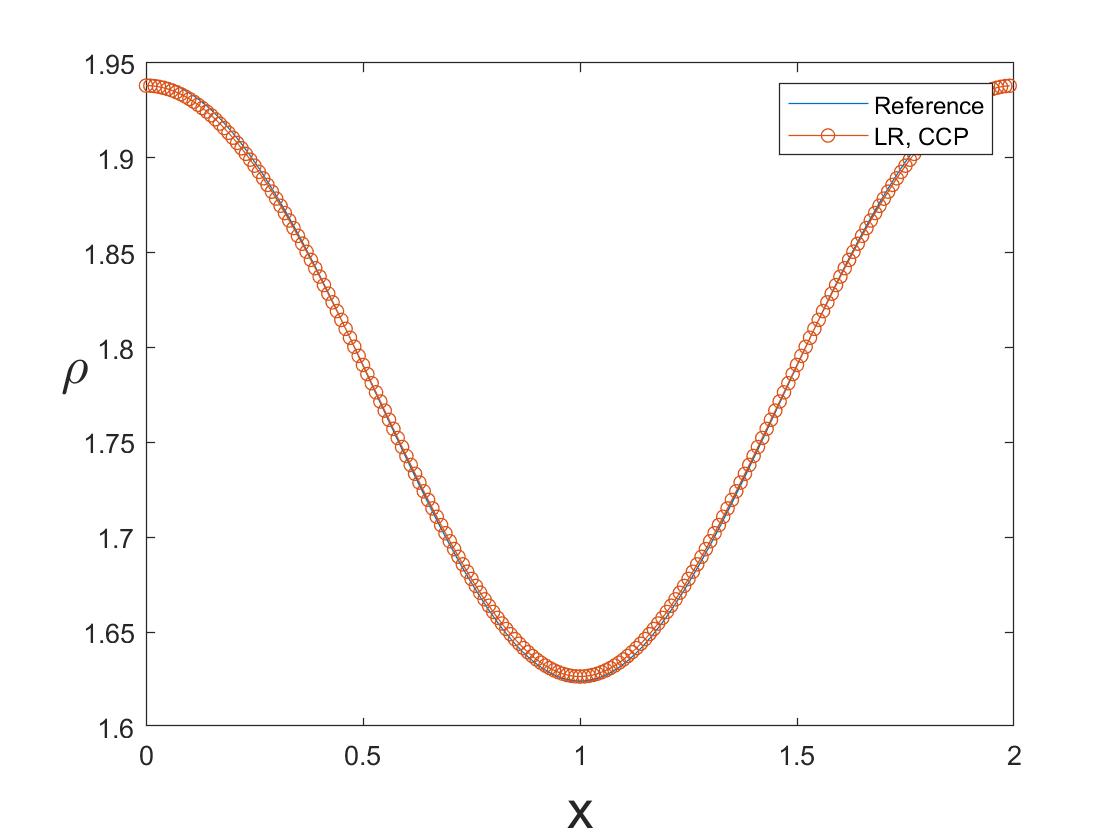}}
     \subfloat{\includegraphics[width=0.47\textwidth]{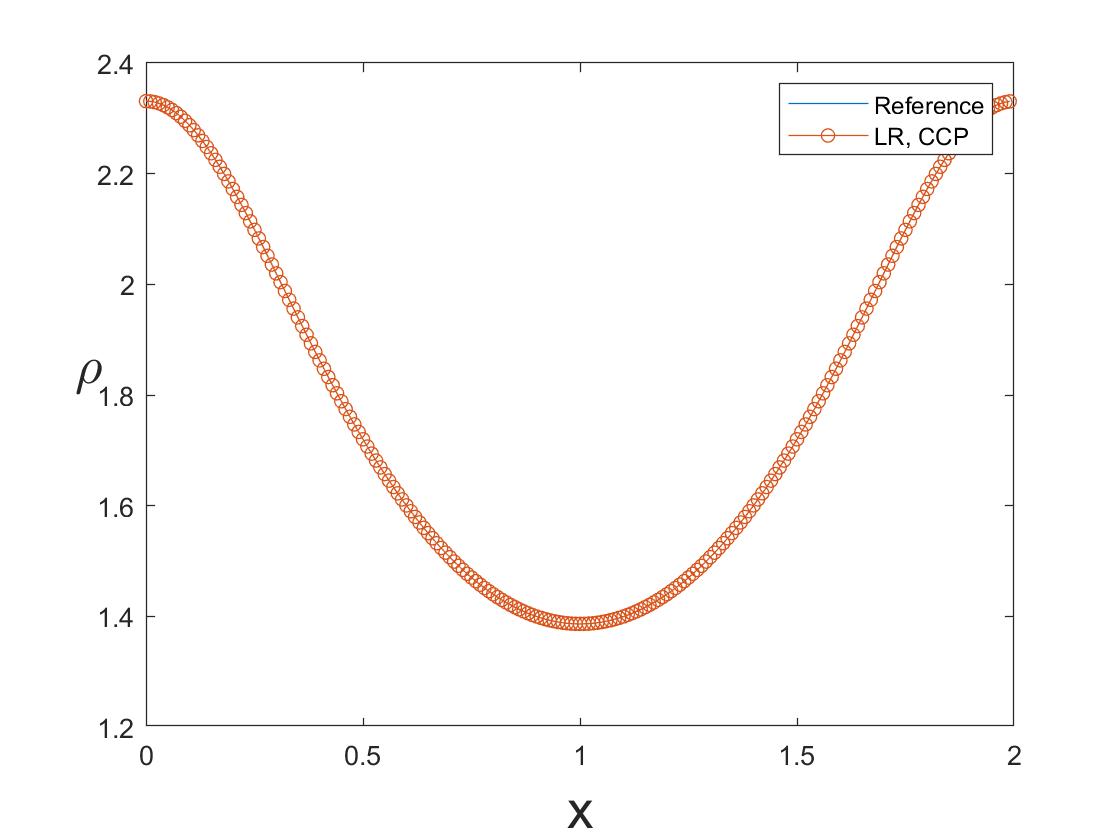}}
     \caption{Example $\mathrm{I}$. We set \revise{$N_x=200, N_v=100, \Delta t=\Delta t_2=\Delta x/3, \Delta t_1=(\Delta t_2)^2, r=20$}. Left: $\epsilon=1$, we compare the density $\rho$ from the dynamical low-rank algorithm with a reference solution at $t_\text{max}=1$. Right: $\epsilon=10^{-3}$, we compare the density $\rho$ from the dynamical low-rank algorithm with diffusion limit at $t_\text{max}=0.1$.}
     \label{Figure1}
\end{figure}

We consider a toy problem with a constant scattering cross section
\[
\sigma(x)=2,\quad x\in [0,2]
\]
and the initial condition set as
\begin{equation}\label{Intialcondition1}
f(0,x,v)=\left((x-1)^2+1\right)\left(v^2+1\right).
\end{equation}
In running Algorithm \ref{Algorithm1}, the dynamical low-rank approximation, we combines central difference and upwind (CCP) flux for $\mathsf{D}$, namely:
\begin{equation}\label{CCP}
  \left(\mathsf{D}u^n\right)_{i,j}=\left\{
  \begin{aligned}
  \epsilon\frac{u^n_{i}-u^n_{i-1}}{\Delta x}+(1-\epsilon)\frac{u^n_{i+1}-u^n_{i-1}}{2\Delta x},\ v_j>0\,,\\
  \epsilon\frac{u^n_{i+1}-u^n_{i}}{\Delta x}+(1-\epsilon)\frac{u^n_{i+1}-u^n_{i-1}}{2\Delta x},\ v_j\leq 0\,.\\
  \end{aligned}\right.
\end{equation}
The combination factor is determined by the Knudsen number $\epsilon$. This means in the kinetic regime when $\epsilon$ is close to $1$ the flux becomes purely upwind but in the diffusion regime when $\epsilon\to 0$, the scheme is of central type. In the kinetic regime for $\epsilon = 1$, we set \revise{$N_x=200$, $N_v=100$, and $\Delta t=\Delta x/3$} for computing the reference solution, and we compare the low rank integrator solution with $r=20$ and the same $(N_x,N_v)$ to this reference solution. \revise{Furthermore, $\Delta t_2 = \Delta x/3$ and $\Delta t_1 = (\Delta t_2)^2$.} In the diffusion regime, for $\epsilon=10^{-3}$, the reference solution is given by the numerical solution to the diffusion equation directly. Both cases are shown in Figure \ref{Figure1} (for $t_\text{max}=1$ and $t_\text{max}=0.1$ respectively) where we clearly see that the numerical solution matches the reference.

\begin{figure}
     \centering
     \subfloat{\includegraphics[width=0.47\textwidth]{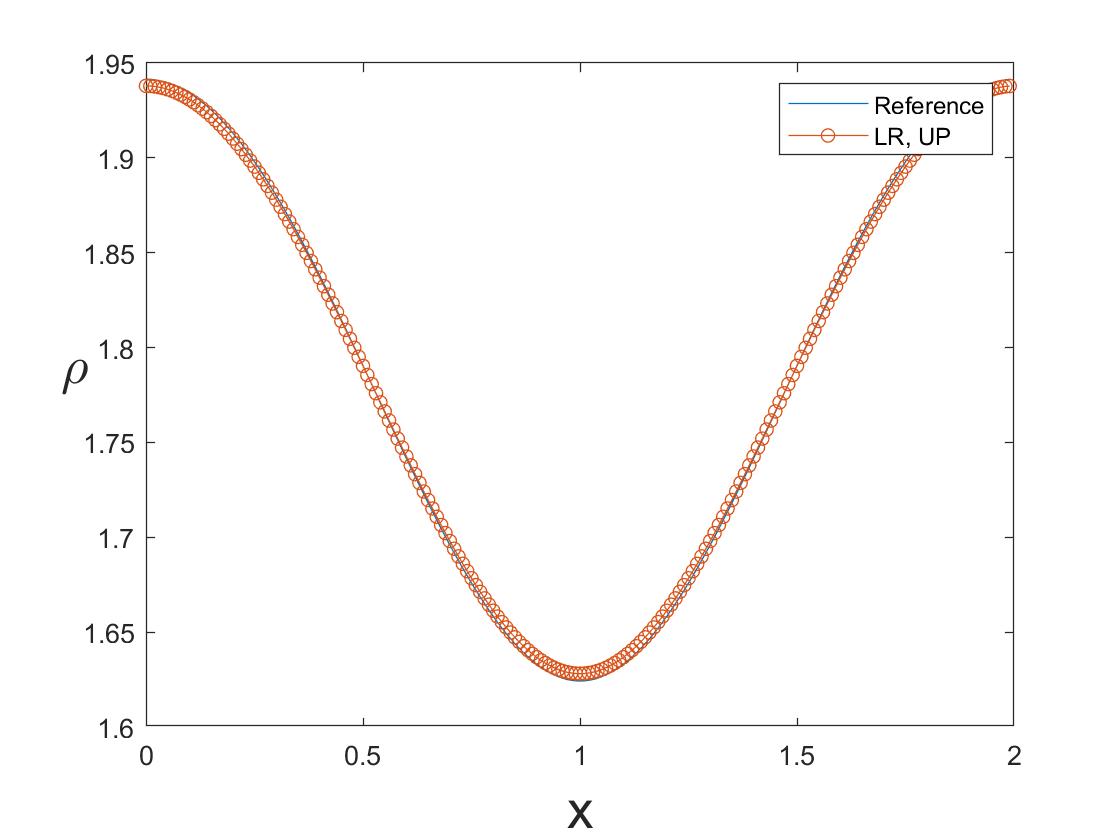}}
     \subfloat{\includegraphics[width=0.47\textwidth]{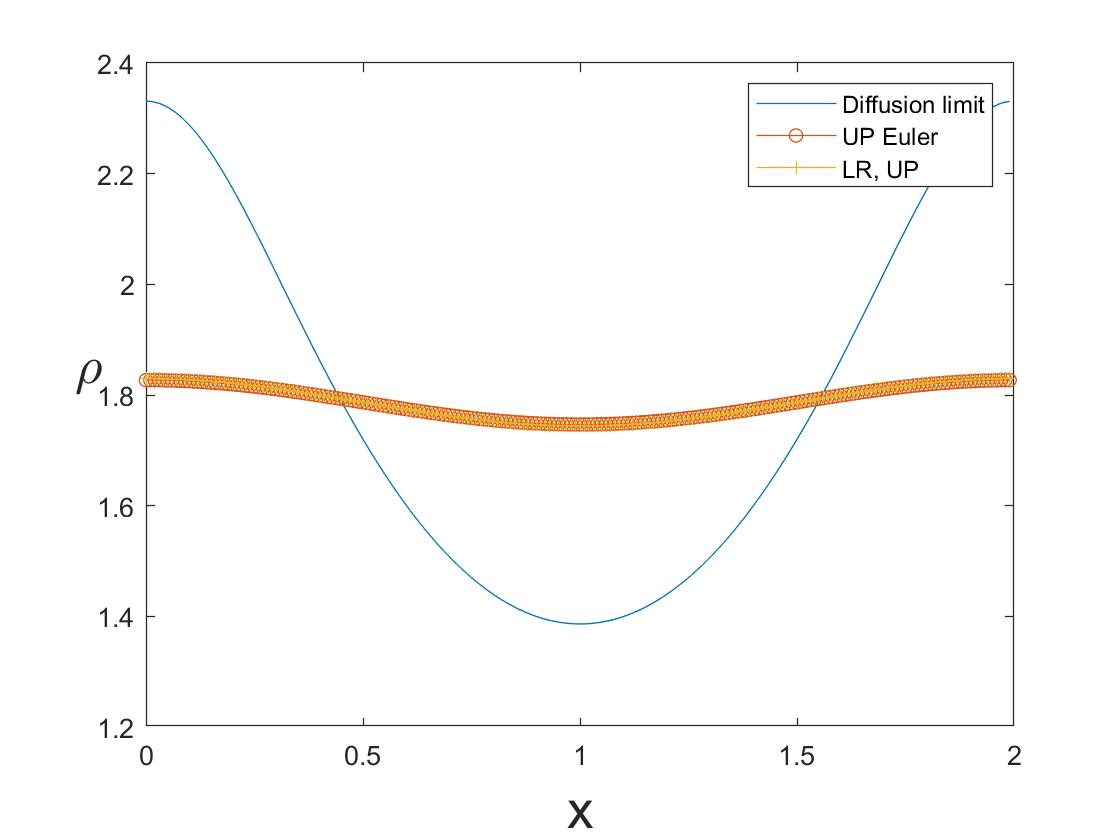}}
     \caption{Example $\mathrm{I}$. \revise{$N_x=200, N_v=100, \Delta t=\Delta t_2=\Delta x/3, \Delta t_1=(\Delta t_2)^2, r=20$}. Left: $\epsilon=1$, we compare the reference solution and the dynamical low-rank integrator solution with upwinding type of flux at $t_\text{max}=1$. Right: $\epsilon=10^{-3}$, we compare the solution to the diffusion equation with the dynamical low-rank integrator solution with upwinding type of flux at $t_\text{max}=0.1$. The $\mathsf{D}(\MSSigma^{-1}\mathsf{D})$ term loses its symmetry, and the reference solution is not captured.}
     \label{Figure2}
\end{figure}

We realize the upwind flux typically brings high artificial diffusion, and this diffusion would be magnified in the diffusion equation when the flux term becomes stiff. To justify the choice of flux $\mathsf{D}$ defined above, here we compute the solution using a dynamical low-rank approximation with $\mathsf{D}$ simply set as the upwind type. The results are shown in Figure \ref{Figure2}. For relatively big Knudsen number ($\epsilon=1$), the low rank integrator solution still agrees with the reference solution well, but the behavior significantly deteriorates in the diffusion regime when $\epsilon\to0$. This is expected as stated in Remark~\ref{rmk:choice_D} that the stencil for the upwind scheme is not symmetric, leading to the fact that $\mathsf{D}(\MSSigma^{-1}\mathsf{D})$ is not a self-adjoint operator as it should be for the $\epsilon\to0$ limit.

\subsection{Example $\mathrm{II}$} 
\begin{figure}
     \centering
     \subfloat{\includegraphics[width=0.47\textwidth]{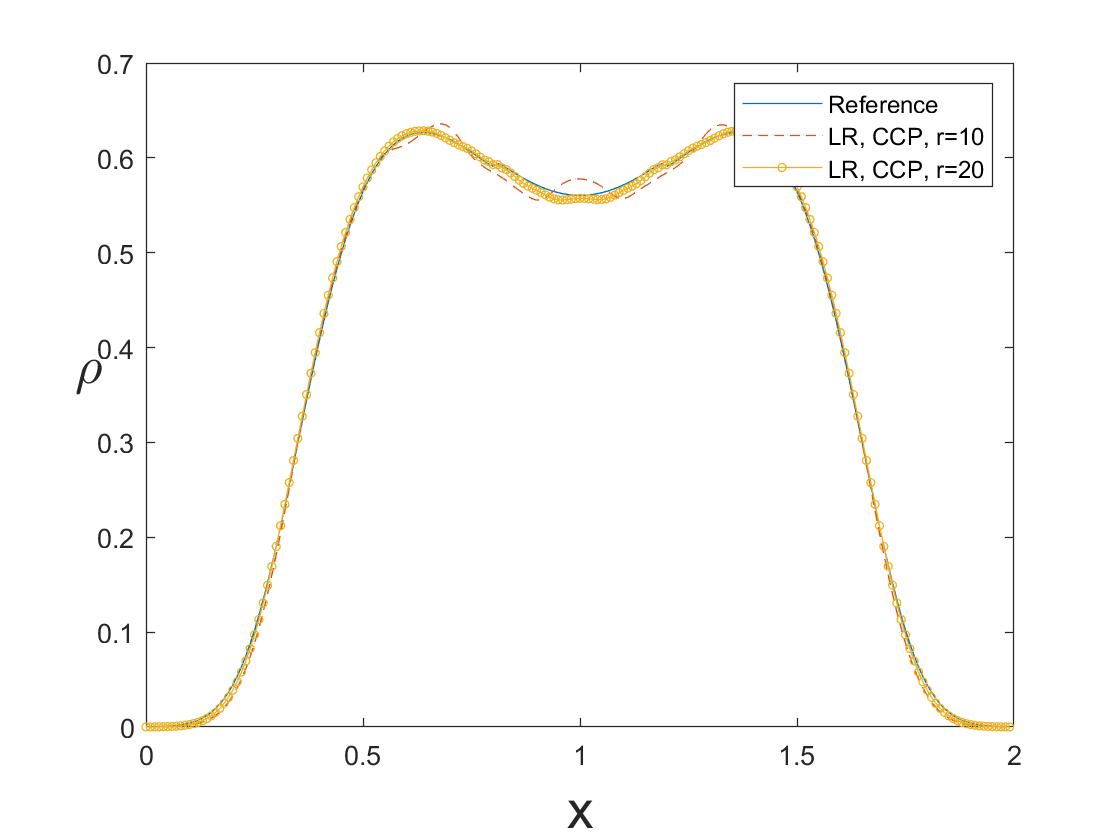}}
     \subfloat{\includegraphics[width=0.47\textwidth]{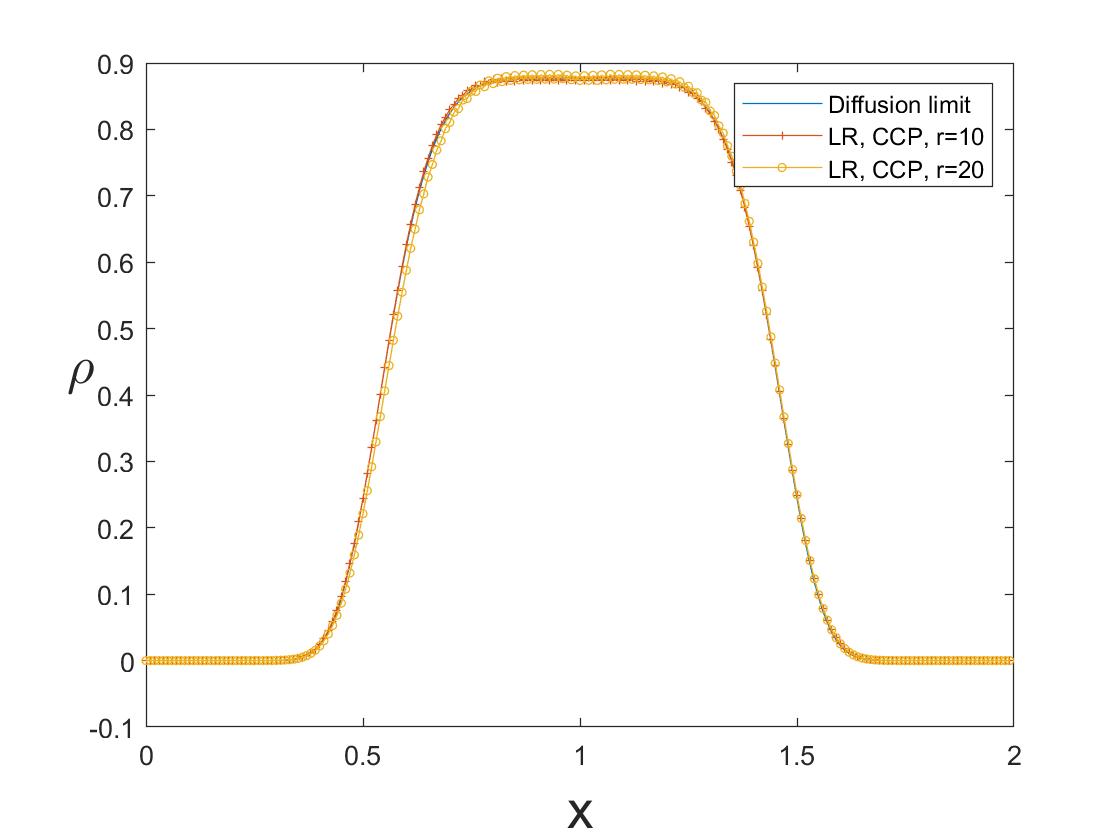}}
     \caption{Example $\mathrm{II}$. \revise{$N_x=200, N_v=100, \Delta t=\Delta t_2=\Delta x/3, \Delta t_1=(\Delta t_2)^2$}. Left: $\epsilon=1$ and $r=10,20$, we compare the density $\rho$ from the low-rank algorithm using CCP for $\MSD_k$ with the implicit Euler/upwind solver at $t_\text{max}=1$. Right: $\epsilon=10^{-3}$ and $r=10,20$, we compare the density $\rho$ from the low-rank algorithm using CCP with the diffusion limit at $t_\text{max}=0.1$. In the plots, LR stands for dynamical low-rank solution.}
     \label{Figure3}
\end{figure}

\begin{figure}
     \centering
     \subfloat{\includegraphics[width=0.47\textwidth]{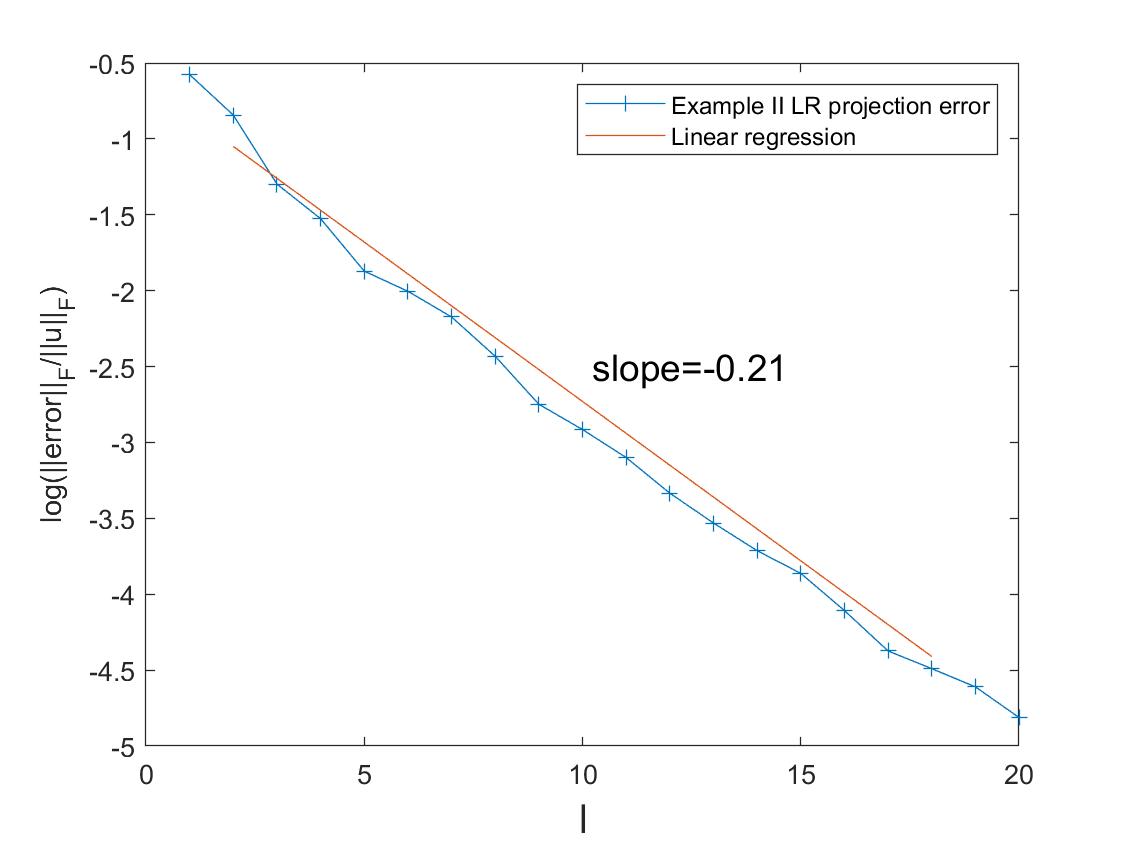}}
     \subfloat{\includegraphics[width=0.47\textwidth]{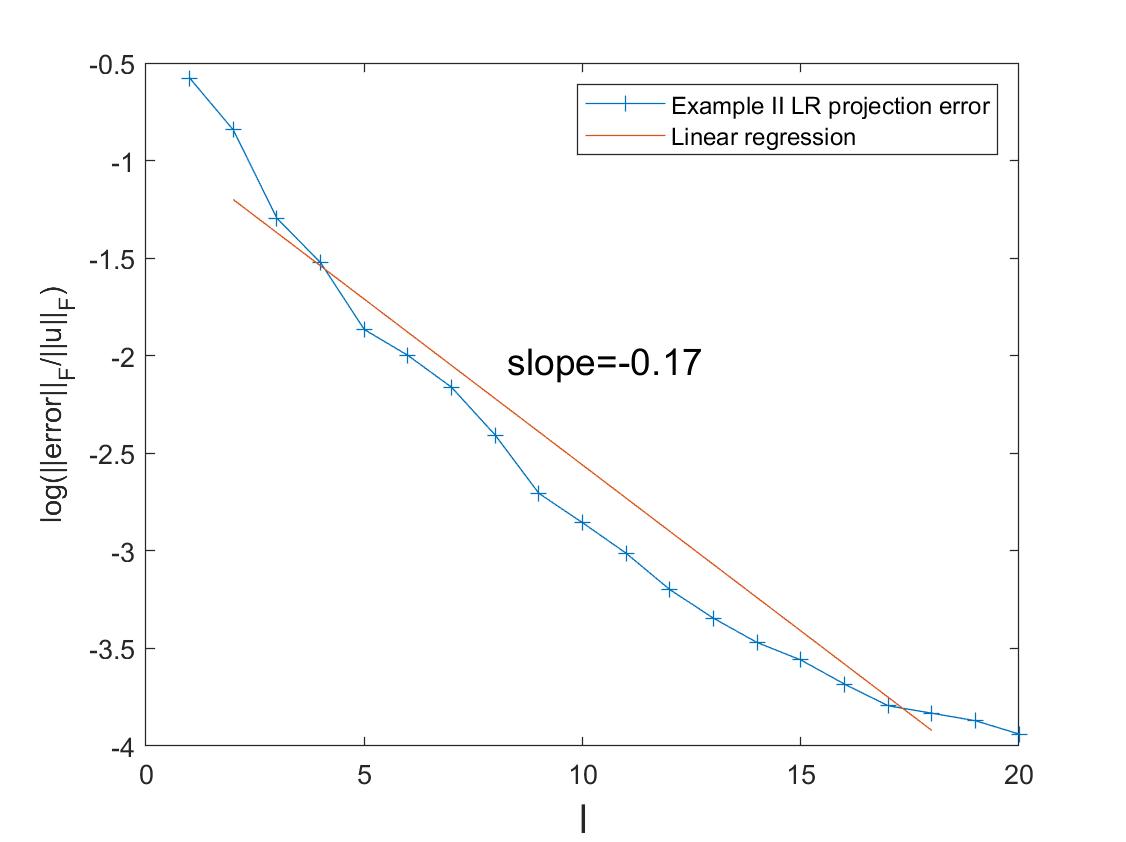}}
     \caption{Example $\mathrm{II}$. \revise{$N_x=200, N_v=100, \Delta t=\Delta t_2=\Delta x/3, \Delta t_1=(\Delta t_2)^2, r=20$.} The numerical solution $u$ is computed using the implicit Euler/upwind solver and $\mathsf{X},\mathsf{V}$ is computed using the low-rank algorithm. Left: $\log\left(\frac{\mathcal{P}^{error}_{x,l}}{\|u\|_F}\right)$ as a function of $l$. Right: $\log\left(\frac{\mathcal{P}^{error}_{v,l}}{\|u\|_F}\right)$ as a function of $l$.}
     \label{Figure8}
\end{figure}

\begin{figure}
     \centering
     \subfloat{\includegraphics[width=0.47\textwidth]{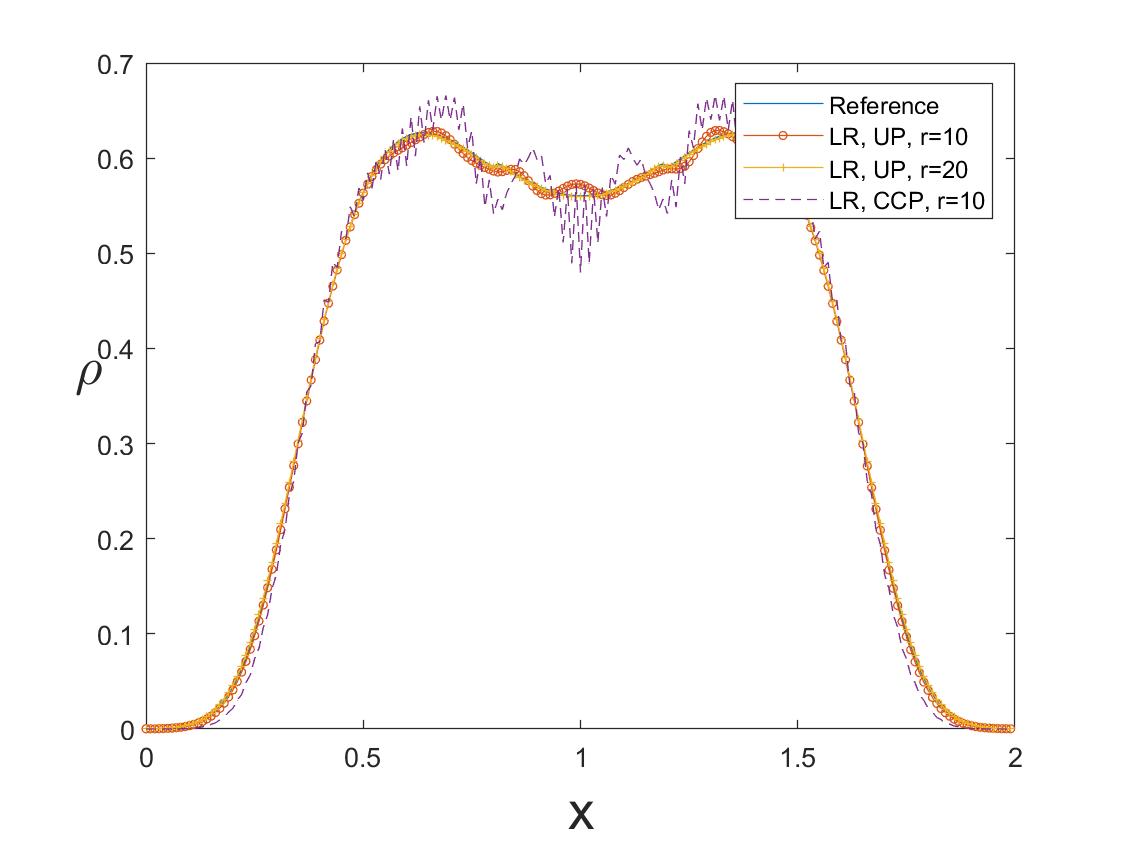}}
     \subfloat{\includegraphics[width=0.47\textwidth]{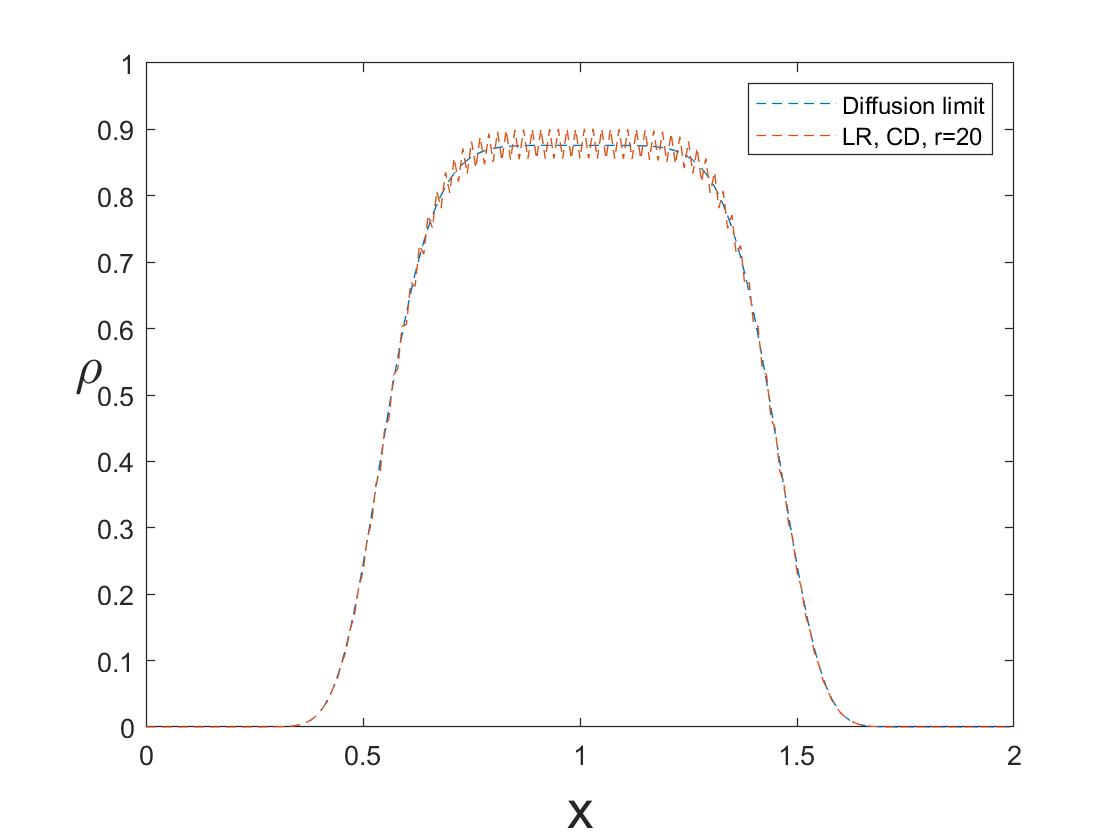}}
     \caption{Example $\mathrm{II}$. \revise{$N_x=200, N_v=100, \Delta t=\Delta t_2=\Delta x/3, \Delta t_1=(\Delta t_2)^2$}. Left: $\epsilon=1$ and $r=10,20$, we compare the density $\rho$ from the low-rank algorithm using upwind and central differences for $\MSD_k$ with the implicit Euler/upwind solver at $t_\text{max}=1$. Right: $\epsilon=10^{-3}$ and $r=20$, we compare the density $\rho$ from the low-rank algorithm using central differences with the diffusion limit at $t_\text{max}=0.1$. In the plots, LR stands for dynamical low-rank solution, and CS stands for central scheme, while UP stands for upwinding.}
     \label{Figure4}
\end{figure}

We consider the following initial condition:
\begin{equation}\label{Intialcondition2}
f(0,x,v)=\left\{
	\begin{aligned}
     &\ 2,\quad 0.8<x<1.2\\
     &\ 0,\quad \text{otherwise},
	\end{aligned}
	 \right.
\end{equation}
and isotropic scattering with a cross section
\begin{equation}\label{eqn:sigma_example_2}
    \sigma(x)=100(x-1)^4\,.
\end{equation}
The cross section thus becomes critical at $x=1$. Comparing the dynamical low-rank solution using $r=20$ to the reference solution (fine discretization for $\epsilon=1$ and diffusion limit for $\epsilon = 10^{-3}$), we see good agreement, as shown in Figure \ref{Figure3}. To be quantitative, we also plot the error, defined in~\eqref{eqn:error_def_num}, in Figure~\ref{Figure8}. It is clear that the error decays exponentially fast. We also test the dynamical low-rank integrator with $\mathsf{D}$ set to be of central scheme type. In this case, however, the method provides lots of artificial oscillation, as shown in Figure~\ref{Figure4}. These artificial oscillations do seem to capture the reference solution weakly (with oscillations centered around the true solution).

\subsection{Example $\mathrm{III}$} In the third example, we use the same initial condition as in Example $\mathrm{II}$, but modify the cross section:
\begin{equation}\label{eqn:sigma_example_3}
    \sigma(x)=\left\{
    \begin{aligned}
    &\ 0.02,\quad x\in[0.35,0.65]\cup[1.35,1.65]\,,\\
    &\ 1,\quad x\in[0,0.35)\cup(0.65,1.35)\cup(1.65,2]\,.
    \end{aligned}\right.
\end{equation}

This cross section is of high contrast and has a discontinuity, and thus the equation quickly achieves equilibrium in the optical thick region ($\sigma\sim 1$), while still staying in the kinetic regime in the optical thin region ($\sigma\sim 0.02$). The dynamical low-rank integrator uses $r=20$. To obtain the reference solutions, we either compute the equation with fine grid when $\epsilon = 1$, or compute the diffusion equation directly with $\epsilon = 10^{-3}$. We also compute the error: it decays exponentially fast, as plotted in Figure~\ref{Figure9}.
\begin{figure}
     \centering
     \subfloat{\includegraphics[width=0.47\textwidth]{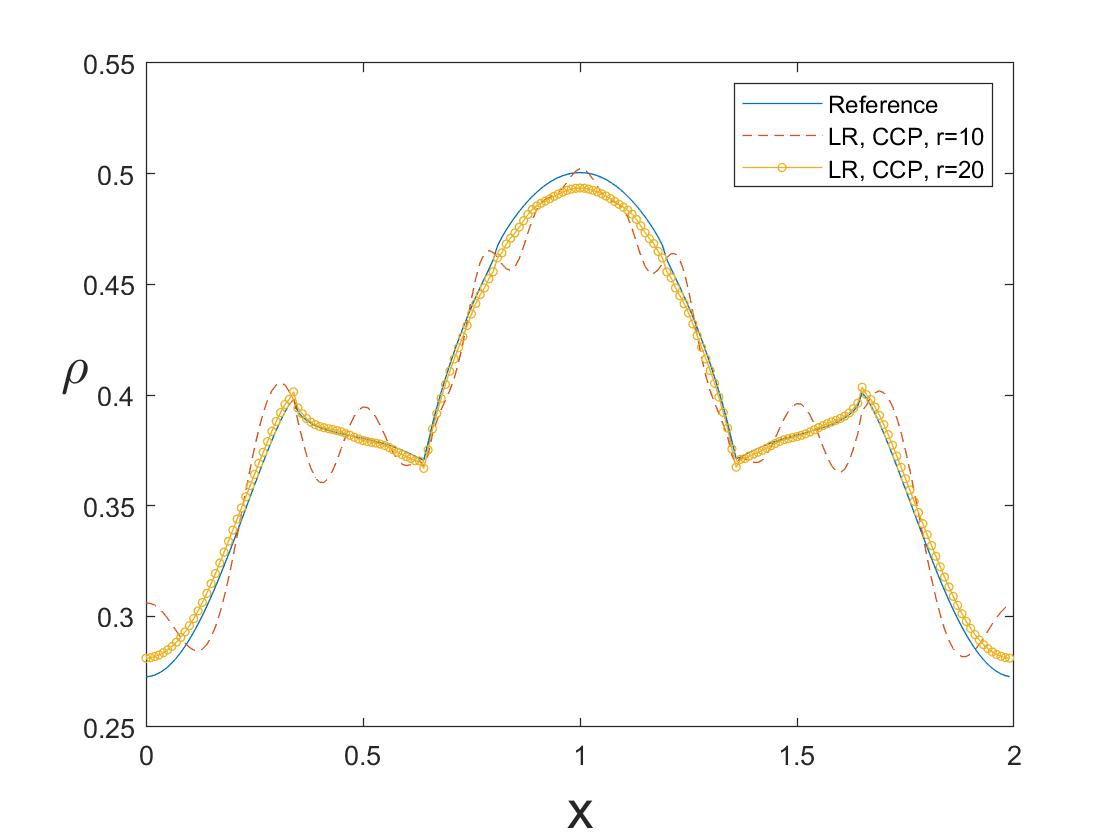}}
     \subfloat{\includegraphics[width=0.47\textwidth]{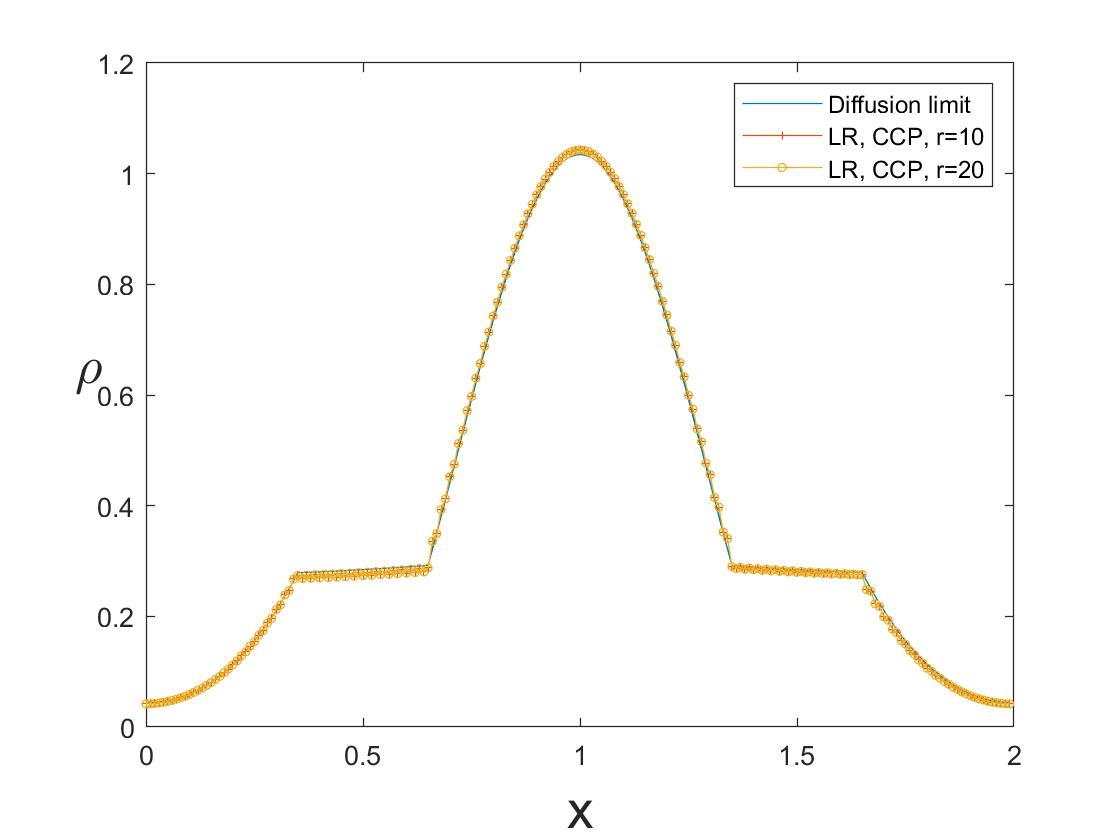}}
     \caption{Example $\mathrm{III}$. \revise{$N_x=200, N_v=100, \Delta t=\Delta t_2=\Delta x/3, \Delta t_1=(\Delta t_2)^2$}. Left: $\epsilon=1$ and $r=10,20$, we compare the density $\rho$ from the low-rank algorithm using CCP for $\MSD_k$ with the implicit Euler/upwind  solver at $t_\text{max}=1$. Right: $\epsilon=10^{-3}$ and $r=10,20$, we compare the density $\rho$ from the low-rank algorithm using CCP with the diffusion limit at $t_\text{max}=0.1$.}
     \label{Figure5}
\end{figure}

\begin{figure}
     \centering
     \subfloat{\includegraphics[width=0.47\textwidth]{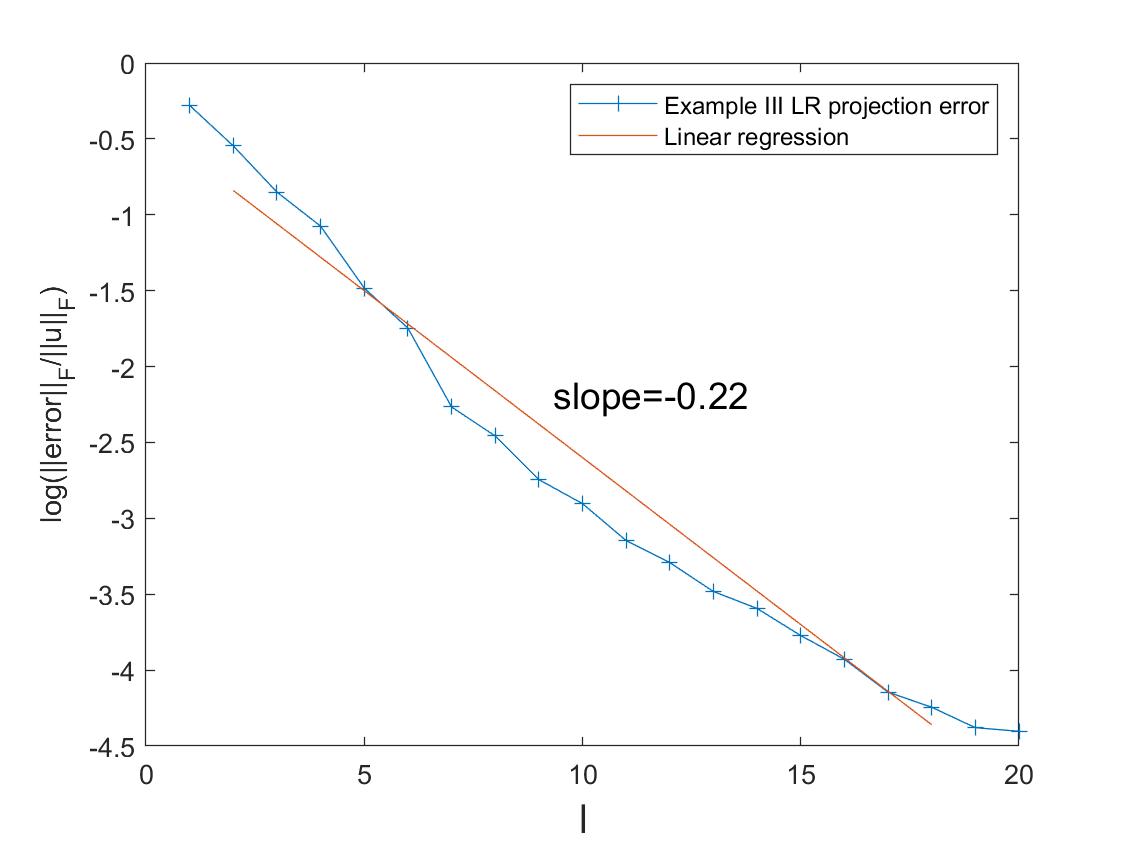}}
     \subfloat{\includegraphics[width=0.47\textwidth]{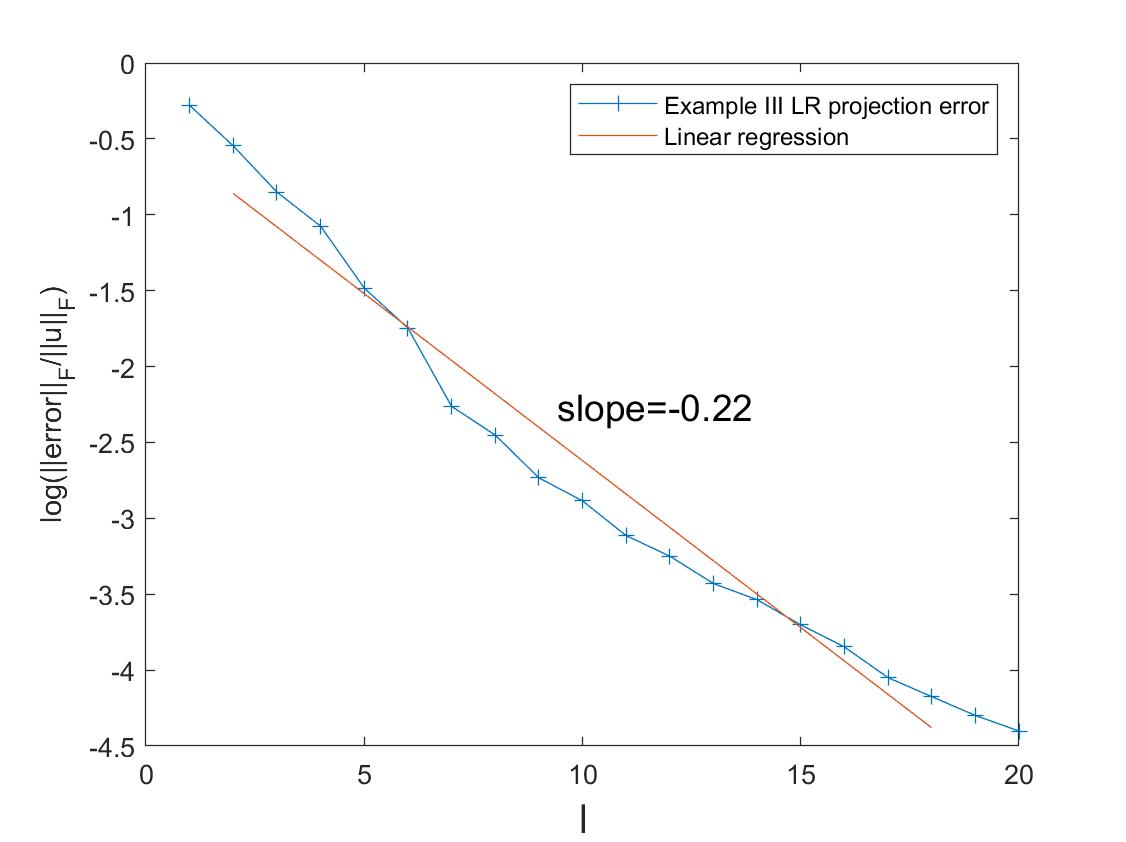}}
     \caption{Example $\mathrm{III}$. Here we use \revise{$N_x=200, N_v=100, \Delta t=\Delta t_2=\Delta x/3, \Delta t_1=(\Delta t_2)^2, r=20$}. the numerical solution $u$ is computed using the implicit Euler/upwind solver and $\mathsf{X},\mathsf{V}$ is computed using the low-rank algorithm. Left: $\log\left(\frac{\mathcal{P}^{error}_{x,l}}{\|u\|_F}\right)$ as a function of $l$. Right: $\log\left(\frac{\mathcal{P}^{error}_{v,l}}{\|u\|_F}\right)$ as a function of  $l$.}
     \label{Figure9}
\end{figure}

\section{Conclusion}
We have applied a projector splitting based dynamical low-rank approximation to a multi-scale linear Boltzmann equation in diffusive scaling. Theoretically the equation has a low rank structure in this regime, and this paper studies if the dynamical low-rank numerical method, a highly nonlinear numerical method, can capture such structure.

The result crucially depends on the employed time integrator. As shown above, if the implicit Euler method is employed, the time stepsize needs to be small. However, the \revise{CNIE} method, due to its embedded symmetry, can preserve the low rank structure automatically for mild time stepsize requirements, if the initial data is well-prepared. Combining these conclusions, we propose to run the Euler method for a very small time stepsize for merely one step to enforce the low rank structure, and then switch to \revise{CNIE} to preserve the rank structure. Numerical evidence is shown to agrees well with the theoretical results obtained.

To our knowledge, this is the first analytical result on showing that the project-splitting dynamical low-rank approximation method truly preserves the rank structure in the kinetic framework.

%An asymptotic analysis has been performed to investigate under which conditions this numerical scheme is preserves the low rank structure of the equation. This crucially depends on the time integrator employed in solving the equations for the low-rank factors. The implicit Euler scheme is AP only if the time step size is very small. The Crank-Nicolson method, on the contrary, does not have such restriction for preserving the asymptotic limit when the initial data is well-prepared, but is not able to drive the solution from a non-equilibrium state to the appropriate low-rank space in time. Thus, we did suggest a combination of both methods in order to obtain an AP method that is independent of the step size and still automatically captures the correct low-rank limit. We have then made this analysis rigorous, which requires slightly more stringent conditions on the initial data. To our knowledge, these are the first AP results for a projector splitting based dynamical low-rank algorithm. Numerical simulation have been conducted that agree with the theoretical results obtained.

\vskip1cm
\noindent\textbf{Acknowledgement.}
The work of QL is supported in part by National Science Foundation under the grant DMS-1619778, and RNMS KI-NET 1107291. The work of ZD is supported in part by Wisconsin Data Science Initiative, and National Science Foundation under the grant DMS-1750488 and TRIPODS: 1740707.

\bibliographystyle{plain}
\bibliography{low_rank_RTE.bib}

\newpage
\appendix
\section{Details for the proof of Theorem \ref{Theorem2}}\label{Theorem2detail}
\begin{lemma}\label{lemmanew}
If Assumption \ref{Assump0} holds true, then for all $n>0$, $\MSL^n_0$ and $\MSS^{n+1/3}_0$ can be written as follows
\begin{equation}\label{AFforLS}
\MSL^n_0=l^n_0e^{\top}_\n,\quad \MSS^{n+1/3}_0=s^{n+1/3}_0\left(\alpha^{n+1}\right)^{\top},\quad \alpha^{n+1}=\MSV^{n+1}e^{\top}_\n,
\end{equation}
with some $l^n_0\in\mathbb{R}^{r\times1}$ and $s^{n+1/3}_0\in\mathbb{R}^{r\times1}$.
\end{lemma}
\begin{remark} Using Assumption \ref{Assump0} and Lemma \ref{lemmanew}, we can further prove that in each step, $\MSL^{n+1/3}$ can be written as
\begin{equation}\label{AMSLnew}
\MSL^{n+1/3}=\hat{l}^{n+1/3}e^{\top}_\n-\epsilon\sum^d_{k=1}(A^n_{\sigma})^{-1}A^{n}_{\partial_k}\hat{l}^{n+1/3}e^{\top}_\n\MSPi_{v_k}+O(\epsilon^2),
\end{equation}
    where $\hat{l}^{n+1/3}=l^{n+1/3}_0+\epsilon l^{n+1/3}_1$. 
\end{remark}
\begin{proof}
The proof proceeds by induction. First, for $n=1$, by \eqref{K1}, after doing QR decomposition, we obtain
\begin{equation}\label{AFirstCP1}
\MSS^{1}_0=s^{1}e_{1,r}(\alpha^1)^{\top},\quad e_{1,r}=\left(1,0,\dots,0\right)^{\top}\in\mathbb{R}^{r\times1},
\end{equation}
where $s^{1}\in\mathbb{R}$ and $\alpha^1=(\MSV^1)^{\top}e_\n$. Using Assumption \ref{Assump0}, we then get
\begin{equation}\label{AFirstCP2}
\MSL^1_0=s^1e_{1,r}(\alpha^1)^{\top}(\MSV^1)^{\top}+O(\epsilon)=s^1e_{1,r}(\alpha^{1,*})^{\top}(\MSV^{1,*})^{\top}+O(\epsilon)=s^1e_{1,r}e^{\top}_\n+O(\epsilon) 
\end{equation}
Let $l^1_0=s^1e_{1,r}$, if we consider terms of order $O(1/\epsilon)$ in the asymptotic expansion of equation \eqref{FDLnew}, we obtain
\begin{equation}
\MSL^{4/3}_0+\MSL^1_0=l^{4/3}_0e^{\top}_\n+l^{1}_0e^{\top}_\n.
\end{equation}
After performing the QR decomposition we have
\begin{equation}\label{ASecondCP}
\MSS^{4/3}_0=\MSL^{n+1/3}_0\MSV^{2}+O(\epsilon)=l^{4/3}_0e^{\top}_\n\MSV^2+O(\epsilon)=l^{4/3}_0\left(\alpha^2\right)^{\top}+O(\epsilon).
\end{equation}
This shows that equation \eqref{AFforLS} holds true for $n=1$.

Now if for $n=k-1$ equation \eqref{AFforLS} holds true, then we can show that
\begin{equation}\label{ALlaststep}
\MSL^{k}_0=l^k_0e^{\top}_\n,
\end{equation}
which implies 
\begin{equation}
\MSL^{k+1/3}_0+\MSL^{k}_0=l^{k+1/3}_0e^{\top}_\n+l^k_0e^{\top}_\n.
\end{equation}
Then we perform a QR decomposition and obtain
\begin{equation}
\MSS^{k+1/3}_0=s^{k+1/3}_0\left(\alpha^{n+1}\right)^{\top},
\end{equation}
where $s^{k+1/3}_0=l^{k+1/3}_0$.
\end{proof}

We now, consider the proof of Theorem \ref{Theorem2} in more detail.
\begin{proof}
\textbf{Step 1-2:}
Instead of \eqref{FinalL} and \eqref{FinalS}. We can obtain
\begin{equation}\label{AnewFinalL}
\frac{u^{n+1/3}_0e-u^n_0e}{\Delta t}-\frac{1}{d}\sum^{d}_{k=1}\MSX^n\MSA^n_{\partial_k}(\MSA^n_{\sigma})^{-1}\MSA^{n}_{\partial_k}(\MSX^n)^{\top}\frac{u^{n+1/3}_0e+u^n_0e}{2}=0,
\end{equation}
\begin{equation}\label{AnewFinalS}
\frac{u^{n+2/3}_0e-u^{n+1/3}_0e}{\Delta t}+\frac{1}{d}\sum^{d}_{k=1}\MSX^n\MSA^n_{\partial_k}(\MSA^n_{\sigma})^{-1}\MSA^{n}_{\partial_k}(\MSX^n)^{\top}\frac{u^{n+2/3}_0e+u^{n+1/3}_0e}{2}=0
\end{equation}
This gives us
\begin{align}
&\rho^{n+1/3}_0=\left(I-\frac{(\Delta t)}{2d}\mathcal{L}^n\right)^{-1}\left(I+\frac{(\Delta t)}{2d}\mathcal{L}^n\right)\rho^n_0,\\
&\rho^{n+2/3}_0=\left(I+\frac{(\Delta t)}{2d}\mathcal{L}^n\right)^{-1}\left(I-\frac{(\Delta t)}{2d}\mathcal{L}^n\right)\rho^{n+1/3}_0\\
\Rightarrow &\rho^{n+2/3}_0=\rho^{n}_0.
\end{align}
Now, we obtain
\begin{equation}
\rho^{n+2/3}_0=\rho^{n}_0,\quad \text{for}\ n>0.
\end{equation}
in the first two splitting step with \revise{CNIE} scheme.

\textbf{Step 3:}
Using equation \eqref{AFforLS} we have
\begin{equation}
\MSS^{n+1/3}_0=s^{n+1/3}_0\left(\alpha^{n+1}\right)^{\top}.
\end{equation}
From the terms of order $O(1/\epsilon)$ we then obtain
\begin{equation}
\MSS^{n+2/3}_0=s^{n+2/3}_0\left(\alpha^{n+1}\right)^{\top}\Rightarrow \MSK^{n+2/3}_0=k^{n+2/3}_0\left(\alpha^{n+1}\right)^{\top}.
\end{equation}
Repeating the procedure used to derive equations \eqref{E1.13}-\eqref{AE1.21}, we still recover the diffusion limit in the last splitting step as
\begin{equation}
\frac{\rho^{n+1}_0-\rho_0^{n+2/3}}{\Delta t}-\frac{1}{d}\sum^d_{k=1}\MSD_{{k}}\left(\MSSigma^{-1}\MSD_{{k}}\rho^{n+1}_0\right)=0.
\end{equation}

\end{proof}

\section{Discussion on the technical assumption}\label{SectionAssump}

We first note that conditions \eqref{Aorth1},\eqref{Aorth2} in Assumption \ref{Assump0} guarantee the orthonormality of $V^{1,*}$. Now, using  equation \eqref{Aorth2} and denoting $\alpha^{n+1,*}=\MSQ^{\top}\alpha^{n+1}$ we have
\begin{equation}\label{AE1.12}
\alpha^{n+1,\ast}(\alpha^{n+1,\ast})^{\top}=\MSQ^{\top}\alpha^{n+1}(\alpha^{n+1})^{\top}\MSQ,\ (\alpha^{n+1,\ast})^{\top}\alpha^{n+1,\ast}=(\alpha^{n+1})^{\top}\alpha^{n+1}=1+O(\epsilon^4).
\end{equation}
For each $1\leq k\leq d$, we define
\[\MSXi^{n+1,\ast}_{v_k}=(\MSV^{n+1,\ast})^{\top}diag(\mathcal{V}_k)\MSV^{n+1,\ast}\,,\]
then we get
\begin{equation}\label{ACVSTAR}
\MSXi^{n+1,\ast}_{v_k}=\MSQ^{\top}\MSXi^{n+1}_{v_k}\MSQ,\quad \left(\MSXi^{n+1,\ast}_{v_k}\right)_{1,j}=\left(\MSXi^{n+1,\ast}_{v_k}\right)_{j,1}=\frac{1}{\sqrt{d}}\delta_{k+1,j}+O(\epsilon),\quad 1\leq j\leq N_v,
\end{equation}
where the second equality can be shown by using equation \eqref{Aorth2}.

By considering equations \eqref{E1.11} and \eqref{AMSLnew}, it suffices to show that the space spanned by the rows of $l^{n+1/3}_0e^{\top}_\n-\epsilon\sum^d_{i=1}(\MSA^n_{\sigma})^{-1}\MSA^{n}_{\partial_i}l^{n+1/3}_0e^{\top}_\n\MSPi_{v_i}$ can contain $e^{\top}_\n,\mathcal{V}^{\top}_1,\dots,\mathcal{V}^{\top}_d$. This, in turn, is equivalent to show that the following $\mathbb{R}^{r\times (d+1)}$ matrix has at least rank $d+1$
\begin{equation}\label{Gamma}
\mathcal{R}=\left(\MSA^n_{\sigma}l^{n+1/3}_0,\MSA^{n}_{\partial_1}l^{n+1/3}_0,\cdots,\MSA^{n}_{\partial_d}l^{n+1/3}_0\right).
\end{equation}
For convenience, we only consider $\sigma(x)=1$ and $\MSD_i$ is skew-symmetric.

First, we define the operator $\mathcal
{H}:\mathbb{R}^{N_x\times1}\rightarrow\mathbb{R}^{d\times d}$ as follows
\begin{equation}\label{Hdefine}
\left(\mathcal{H}(\textbf{k})\right)_{i,j}=\left(\MSD_i\textbf{k}\right)^{\top}\left(\MSD_j\textbf{k}\right),\quad \textbf{k}\in\mathbb{R}^{N_x\times1}
\end{equation}
\begin{lemma}\label{LemmaforAS}
Assume Assumption \ref{Assump0} is true for the numerical solution at time $t_n$, $\frac{\Delta t}{(\Delta x)^2}$ and $\epsilon$ are small enough, if 
\begin{equation}\label{det}
det(\mathcal{H}(k^{n}))\neq0,
\end{equation}
where $k^n$ is defined in equation \eqref{K1}, then Assumption \ref{Assump0} is true for the numerical solution at time $t_{n+1}$.
\end{lemma} 

\begin{remark} By Lemma \ref{LemmaforAS}, we only need to check condition \eqref{det}. Since $k^n$ is a discrete approximation to $\rho(x,t^n)$, assume numerical error is small and the solution satisfies
\[
D(\rho)(t)=det\left(\int_{\Omega_x}\nabla \rho(x,t)\left(\nabla \rho(x,t)\right)^{\top}dx\right)\neq 0,\quad \forall 0\leq t\leq T,
\]
then we have $\mathcal{H}(k^{n})\approx D(\rho)(t^n)$ and \eqref{det}.
\end{remark}
\begin{proof}
If Assumption \ref{Assump0} is true at time $t_{n-1}$, we have
\[
\MSK^{n}=k^{n}\left(\alpha^{n}\right)^{\top}-\epsilon\sum^{d}_{i=1}\MSD_ik^{n}\left(\alpha^{n}\right)^{\top}\MSXi^{n}_{v_i}+O(\epsilon^2).
\]
and there exists an orthogonal matrix $\MSQ$ such that
\begin{equation}\label{semiK}
\MSK^{n}\MSQ=k^{n}\left(\alpha^{n,*}\right)^{\top}-\epsilon\sum^{d}_{i=1}\MSD_ik^{n}\left(\alpha^{n,*}\right)^{\top}\MSXi^{n,*}_{v_i}+O(\epsilon^2),
\end{equation}
where $\alpha^{n,*}$, $\MSXi^{n,*}_{v_i}$ satisfy \eqref{AE1.12} and \eqref{ACVSTAR}.

Use \eqref{AE1.12} and \eqref{ACVSTAR}, we can rewrite \eqref{semiK} as
\begin{equation}\label{semiK2}
\MSK^{n}\MSQ=\left(k^{n},-\frac{\epsilon}{\sqrt{d}}\MSD_1k^{n},\dots,-\frac{\epsilon}{\sqrt{d}}\MSD_dk^{n},\dots,\right)+O(\epsilon^2).
\end{equation}
We notice that $\MSX^{n}$ comes from QR decomposition of $\MSK^{n-1}$. Therefore, use \eqref{semiK2}, there exits an invertible matrix $\MSP\in\mathbb{R}^{r\times r}$ such that 
\begin{equation}\label{newX}
\MSX^{*}=\MSX^{n}\MSP^{-1}=\left(k^{n},-\MSD_1k^{n},\dots,-\MSD_dk^{n},\dots,\right)+O(\epsilon).
\end{equation}
For any $1\leq i\leq d$, we further define
\begin{equation}\label{newAK}
\MSA^*_{\partial_i}=\left(\MSX^{*}\right)^{\top}\MSD_i\MSX^{*},\ \text{then}\ \MSA^n_{\partial_i}=\MSP^{\top}\MSA^*_{\partial_i}\MSP
\end{equation}
and
\begin{equation}\label{NewAsigma}
\MSA^*_{\sigma}=\left(\MSX^{*}\right)^{\top}\MSSigma\MSX^{*},\ \text{then}\ \MSA^n_{\sigma}=\MSP^{\top}\MSA^*_{\sigma}\MSP.
\end{equation}
Plug equations \eqref{newAK},\eqref{NewAsigma} into \eqref{Gamma} and define $l^*=\MSP l^{n+1/3}_0$, it suffices to show that
\[
\mathcal{R}^*=\left(\MSA^*_{\sigma}l^*,\MSA^{*}_{\partial_1}l^*,\cdots,\MSA^*_{\partial_d}l^*\right)
\]
has rank $d+1$. Now we divide the following proof into two steps:

\textbf{First step: \emph{$(l^*)_1\gg 0$.}}  Similar to equation \eqref{ALlaststep}, we have 
\[
\MSS^n\left(\MSV^n\right)^{\top}=\MSL^{n}=l^n_0e^{\top}_\n+O(\epsilon).
\]
Using \eqref{AE1.111},\eqref{K1},\eqref{AE1.12}, we get
\[
\MSX^*\MSP \MSS^n\left(\MSV^n\right)^{\top}=\MSK^{n}\left(\MSV^n\right)^{\top}=k^{n}\left(\alpha^n\right)^{\top}\left(\MSV^n\right)^{\top}=k^{n}\left(\alpha^{n,*}\right)^{\top}\left(\MSV^{n,*}\right)^{\top}=k^ne^{\top}_\n+O(\epsilon^2),
\]
which implies
\[
\MSX^*\MSP l^n_0e^{\top}_\n=k^ne^{\top}_\n.
\]
Because $\MSX^*$ is an invertible matrix and the first column of $\MSX^*$ is $k^n$, we must have 
\begin{equation}\label{ln0}
\left(\MSP l^n_0\right)_1=1+O(\epsilon).
\end{equation}
Recall that
\begin{equation}\label{Lstar0}
\left(l^*\right)_1=\left(\MSP l^{n+1/3}_0\right)_1,
\end{equation}
similar to asymptotic analysis performed to obtain \eqref{FinalL}, we have
\[
\begin{aligned}
&\frac{\MSL^{n+1/3}_0e_\n-\MSL^n_0e_\n}{\Delta t}-\sum^{d}_{i=1}\MSA^n_{\partial_i}(\MSA^n_{\sigma})^{-1}\MSA^{n}_{\partial_i}\frac{l^{n+1/3}_0+l^{n}_0}{2}e^{\top}_\n\MSPi^2_{v_i}e_\n=0\\
\Rightarrow &\frac{l^{n+1/3}_0-l^n_0}{\Delta t}-\frac{1}{d}\sum^{d}_{i=1}\MSA^n_{\partial_i}(\MSA^n_{\sigma})^{-1}\MSA^{n}_{\partial_i}\frac{l^{n+1/3}_0+l^{n}_0}{2}=0,
\end{aligned}
\]
which implies
\begin{equation}\label{Ldifference}
\|l^{n+1/3}_0-l^n_0\|_2=O\left(\frac{\Delta t}{(\Delta x)^2}\right).
\end{equation}
Combining \eqref{ln0},\eqref{Lstar0} with \eqref{Ldifference}, we finally obtain
\[
\left|\left(\MSP\left(l^{n+1/3}_0-l^n_0\right)\right)_1\right|=O\left(\frac{\Delta t}{(\Delta x)^2}\right)+O(\epsilon)\Rightarrow \left(l^*\right)_1>1-\left[O\left(\frac{\Delta t}{(\Delta x)^2}\right)+O(\epsilon)\right]\gg0,
\]
where $\frac{\Delta t}{(\Delta x)^2}$, $\epsilon$ have to be very small.

\textbf{Second step: $\Gamma^*$ is invertible.} Now we can prove \eqref{det} by contradiction. If we assume there exists $\gamma\in\mathbb{R}^{(d+1)\times 1}$ with $\|\gamma\|_2=1$ such that
\begin{equation}\label{zero1}
\mathcal{R}^* \gamma=0\Rightarrow \left(\MSA^*_{\sigma}\gamma_1+
\sum^{d}_{i=i}\MSA^{*}_{\partial_i}\gamma_{i+1}\right)l^*=0.
\end{equation}
 Because $\MSD_i$ is skew-symmetric, $\MSA^{*}_{\partial_i}$ is skew-symmetric for each $i$. From equation \eqref{zero1}, we deduce
\[
(l^*)^{\top}\left(\MSA^*_{\sigma}\gamma_1+
\sum^{d}_{i=i}\MSA^{*}_{\partial_i}\gamma_{i+1}\right)l^*=0\Rightarrow \gamma_1(l^*)^{\top}\MSA^*_{\sigma}l^*=0\Rightarrow \gamma_1=0\Rightarrow 
\sum^{d}_{i=i}\MSA^{*}_{\partial_i}l^*\gamma_{i+1}=0,
\]
where we use $(l^*)^{\top}\MSA^{*}_{\partial_i}l^*=0$. 

Plug equation \eqref{newX} into equation \eqref{newAK}, we obtain
\begin{equation}\label{StructureforAstar}
\left(\MSA^*_{\partial_m}\right)_{i,j}=
\left\{
\begin{aligned}
&\revise{\left(\MSD_{i-1}k^{n}\right)^{\top}\MSD_{m}\left(\MSD_{j-1}k^{n}\right)},\quad 2\leq i,j\leq d+1\\
&0,\quad i=j\\
&-\left(\MSD_{i-1}k^{n}\right)^{\top}\left(\MSD_{m}k^{n}\right)+O(\epsilon),\quad 2\leq i\leq d+1,\ j=1
\end{aligned}
\right.
\end{equation}
for any $1\leq m\leq d$. Consider the first column of $\MSA^*_{\partial_i}$ from equation \eqref{StructureforAstar}. Since $\left(l^*\right)_1\gg0$, equation \eqref{zero1} implies
\[
\mathcal{H}(k^{n})\gamma_{(2:d+1)}=0,
\]
where $\mathcal{H}$ is defined in equation \eqref{Hdefine} and $\gamma_{(2:d+1)}\in\mathbb{R}^{d\times1}$ is the cutoff vector to $\gamma$ defined as
\[
\left(\gamma_{(2:d+1)}\right)_{i}=\gamma_{i+1},\quad 1\leq i\leq d\,.
\]
Since $\mathcal{H}(k^{n})$ is an invertible matrix, this finally shows that $\gamma_{(2:d+1)}$ has to be zero, which contradicts to $\|\gamma\|_{2}=1$.
\end{proof}

\end{document}